\magnification=\magstep1
\input amstex
\documentstyle{amsppt}
\catcode`\@=11 \loadmathfont{rsfs}
\def\mycal{\mathfont@\rsfs}
\csname rsfs \endcsname \catcode`\@=\active

\vsize=6.5in

\topmatter 
\title Independence properties \\in subalgebras of ultraproduct II$_1$ factors
\endtitle
\author SORIN POPA \endauthor

\rightheadtext{Independence in subalgebras}

\affil     {\it  University of California, Los Angeles} \endaffil

\address Math.Dept., UCLA, Los Angeles, CA 90095-1555\endaddress
\email popa\@math.ucla.edu\endemail

\thanks Supported in part by NSF Grant DMS-1101718 and a Simons Fellowship \endthanks

\abstract  Let $M_n$ be a sequence of finite factors with $\dim(M_n)\rightarrow \infty$ and denote $\text{\bf M}=\Pi_\omega M_n$ their ultraproduct over a free ultrafilter $\omega$. We prove that  if  $\text{\bf Q}\subset \text{\bf M}$ is either an ultraproduct $\text{\bf Q}=\Pi_\omega Q_n$ of subalgebras   
$Q_n\subset M_n$, with $Q_n \not\prec_{M_n} Q_n'\cap M_n$, $\forall n$, or the centralizer $\text{\bf Q}=B'\cap \text{\bf M}$ of a separable amenable $^*$-subalgebra 
$B\subset \text{\bf M}$, 
then for any separable subspace $X\subset \text{\bf M}\ominus (\text{\bf Q}'\cap \text{\bf M})$, there exists a diffuse abelian von Neumann subalgebra 
in $\text{\bf Q}$ which is {\it free independent} to $X$, relative to $\text{\bf Q}'\cap \text{\bf M}$. Some related independence properties for subalgebras 
in ultraproduct II$_1$ factors are also discussed.\endabstract

\endtopmatter

\document

\heading 0. Introduction \endheading

We continue in this paper the investigation of  independence properties in subalgebras of ultraproduct II$_1$ factors, from  
[P6], [P12]. The main result we prove along these lines is the following:

\proclaim{0.1. Theorem} Let $M_n$ be a sequence of finite factors with $\dim M_n \rightarrow \infty$ and denote by $\text{\bf M}$ the ultraproduct 
$\text{\rm II}_1$ factor $\Pi_\omega M_n$, over a free ultrafilter $\omega$ on $\Bbb N$. Let $\text{\bf Q}\subset \text{\bf M}$ be a von Neumann 
subalgebra satisfying one of the following: 
\vskip .05in
$(a)$ $\text{\bf Q}=\Pi_\omega Q_n$, for some von Neumann subalgebras $Q_n\subset M_n$ 
satisfying the condition $Q_n \not\prec_{M_n} Q_n'\cap M_n$, $\forall n$ $($in the sense of $\text{\rm [P10]})$; 

$(b)$ $\text{\bf Q}=B'\cap \text{\bf M}$, for some separable amenable von Neumann subalgebra $B\subset \text{\bf M}$. 
\vskip .05in
Then given any separable subspace $X\subset \text{\bf M}\ominus (\text{\bf Q}'\cap \text{\bf M})$, there exists 
a diffuse abelian von Neumann subalgebra $A\subset \text{\bf Q}$ 
such that $A$ is free independent to $X$, relative to $\text{\bf Q}'\cap \text{\bf M}$, i.e. $E_{\text{\bf Q}'\cap \text{\bf M}}(x_0\Pi_{i=1}^n a_ix_i)=0$, 
for all $n\geq 1$, $x_0, x_n \in X\cup \{1\}$, $x_i \in X$, $1\leq i \leq n-1$, $a_i \in A\ominus \Bbb C1$, $1\leq i \leq n$. 
\endproclaim

Note that the particular case when $Q_n\subset M_n$ are II$_1$ factors 
with atomic relative commutant, for which one clearly has $Q_n \not\prec_{M_n} Q_n'\cap M_n$, recovers (2.1 in [P6]).   

The conclusion in 0.1 above can alternatively be interpreted as follows: given any separable von Neumann subalgebra $P$ of $\text{\bf M}$ that makes a commuting 
square with $\text{\bf Q}'\cap \text{\bf M}$ (in the sense of 1.2 in [P2]; see Sec. 1.2 below for the definition) 
and we let $B_1=P\cap (\text{\bf Q}'\cap \text{\bf M})$, there exists a separable von Neumann subalgebra $B_0\subset \text{\bf Q}$, 
such that $P\vee B_0 \simeq P*_{B_1} (B_1\overline{\otimes} B_0)$ (amalgamated free product of finite von Neumann algebras over a common subalgebra, see [V], [P4]). Since in the case $(b)$ of 0.1 we have $\text{\bf Q}'\cap \text{\bf M}=B$
(see 2.1 below) and all embeddings into an ultraproduct II$_1$ factor $\text{\bf M}$ of an amenable separable von Neumann algebra $B$ are conjugate by  unitaries  
in $\text{\bf M}$, Theorem 0.1 shows in particular that if two separable finite von Neumann algebras $N_1, N_2$ 
containing copies of $B$ are embeddable into $\text{\bf M}$, then $N_1*_BN_2$ is embeddable into $\text{\bf M}$ as well.  Note that the case 
$B$ atomic of this result already appears in [P6], while the case $B$ arbitrary but with $\text{\bf M}=R^\omega$
was shown in [BDJ].  More precisely, 0.1 implies the following strengthening of these results:  

\proclaim{0.2. Corollary} Let $\text{\bf M}=\Pi_\omega M_n$ be an ultraproduct $\text{\bf II}_1$ factor as in $0.1$. 
Let $N_i \subset \text{\bf M}$ be separable finite von Neumann subalgebras with   
amenable von Neumann subalgebras $B_i\subset N_i$, $i=1,2$, such that  $(B_1, \tau_{|B_1})\simeq (B_2, \tau_{|B_2})$. 
Then there exists a unitary element $u\in  \text{\bf M}$ so that $uB_1u^*=B_2$ and so that, after identifying $B=B_1\simeq B_2$ this way,  
we have $uN_1u^* \vee N_2 \simeq N_1*_B N_2$.  
\endproclaim

To prove Theorem 0.1, we first construct unitaries $u\in \text{\bf Q}$ that are approximately $n$-independent with respect to given finite sets $X\perp \text{\bf Q}'\cap \text{\bf M}$. 
Taking larger and larger $n$, larger and larger finite sets $X$ and better and better approximations, and combining with a diagonalization procedure, 
one can then get unitaries that are free independent to a given countable set, due to the ultraproduct framework.

The approximately independent  unitary $u$  is constructed by patching together incremental pieces of it, while controlling 
the trace of alternating words involving $u$ and a given set $X$. This technique was initiated in [P3], being then fully developed in [P6],  
where it has been used to prove a particular case of $0.1 (a)$. More recently, it has been used in [P12] to establish existence of 
free independence in ultraproducts of maximal abelian $^*$-subalgebras (abbreviated hereafter MASA) $A_n \subset M_n$ that are {\it 
singular}  in the sense of [D1] (i.e., any unitary element in $M_n$ that normalizes $A_n$ must lie in $A_n$),  thus settling the Kadison-Singer 
problem for the corresponding ultrapower inclusion $\text{\bf A}=\Pi_\omega A_n \subset \Pi_\omega M_n=\text{\bf M}$. 

If in turn the normalizers of the MASAs 
$A_n\subset M_n$ are large, then 
one can still detect certain independence properties inside $\text{\bf A}$, by using the same type of techniques. Thus, 
it was shown in [P12] that $3$-independence always occurs in $\text{\bf A}$, and we prove here that given any countable group of unitaries 
$\Gamma$ in $ \text{\bf M}$, that normalizes $\text{\bf A}$ and  
acts freely on it, there exists a diffuse subalgebra $B_0$ in $\text{\bf A}$ such that  any word $\Pi_{i=1}^n u_ib_iu_i^*$ 
with $b_i\in B_0\ominus \Bbb C1$ and distinct $u_i\in \Gamma$, has trace $0$. This actually amounts 
to $B_0$ being the base of a Bernoulli $\Gamma$-action. We in fact prove the following stronger result:  

\proclaim{0.3. Theorem} Let $A_n \subset M_n$ be MASAs in finite factors, as before, and denote 
$\text{\bf A}=\Pi_\omega A_n \subset \Pi_\omega M_n = \text{\bf M}$. 
Assume $\Gamma \subset \text{\bf M}$ 
is a countable group of unitaries normalizing $\text{\bf A}$ and acting freely on it, and let $H\subset \Gamma$ be an amenable subgroup. 
Given any separable abelian von Neumann 
subalgebra $B \subset \text{\bf A}$, there exists a $\Gamma$-invariant subalgebra $A\subset \text{\bf A}$ such that   
$A, B$ are $\tau$-independent and $\Gamma \curvearrowright A$ is isomorphic to the generalized Bernoulli 
action $\Gamma \curvearrowright L^\infty([0,1]^{\Gamma/H})$. 
\endproclaim 

Note that if the above ultraproduct inclusion $\text{\bf A}\subset \text{\bf M}$ comes from a sequence of finite dimensional 
diagonal inclusions $D_n \subset M_{n\times n}(\Bbb C)$, or is of the form $D^\omega \subset R^\omega$, where $D\subset R$ is the 
unique (up to conjugacy by an automorphism, by [CFW]) Cartan subalgebra of the hyperfinite II$_1$ factor, 
then a countable group $\Gamma$ can be embedded into 
the normalizer $\Cal N_{\text{\bf M}}(\text{\bf A})$ of $\text{\bf A}$  in $\text{\bf M}$, in a way that it acts freely on $\text{\bf A}$, iff 
it is {\it sofic} (in the sense of [W]; see the expository paper [Pe] and [Pa]). Thus, with the terminology in [EL], where 
an {\it action} of a sofic group $\Gamma \curvearrowright X$ 
is called {\it sofic} if the inclusion $L^\infty(X)\subset L^\infty(X)\rtimes \Gamma$ admits a commuting square embedding into $\text{\bf A}\subset \text{\bf M}$, 
with $\Gamma$ embedded into $\Cal N_{\text{\bf M}}(\text{\bf A})$,  it follows  
from  0.3 that if $\Gamma \curvearrowright X$ is sofic then 
any product action $\Gamma \curvearrowright X \times Y$, with $\Gamma\curvearrowright Y=[0,1]^I$ a generalized 
Bernoulli action corresponding to the left action of $\Gamma$ on a set $I= \oplus_i \Gamma/H_i$, for 
some countable family of amenable subgroups $H_i \subset \Gamma$, is sofic. This generalizes a result in [EL].  

The paper is organized as follows. In Section 1 we recall some basic facts needed in the paper, such as the {\it local quantization} 
lemma from [P1], [P5] and the criterion for (non-)conjugacy of subalgebras  from [P10]. We also prove a general fact about 
centralizers (or commutants) of countable sets in ultraproduct II$_1$ factors (see Theorem 1.7). In Section 2 we prove some bicentralizer 
results concerning amenable algebras and groups, in ultrapower framework, that we need in the proofs of 0.1 and respectively 0.3. 
We conjecture that, in fact, the bicentralizer property  characterizes amenability (see 2.3.1). 

In Section 3 we prove the main technical result needed 
in the proof of Theorem 0.1, by using incremental patching techniques. This result, stated as Lemma 3.2, actually amounts to an ``approximate version'' 
of the free independence result in 0.1. In Section 4 we derive Theorem 0.1 (in fact a strengthening of it, stated as Theorem 4.3), 
by using Lemma 3.2 and an appropriate diagonalization procedure. 

In Section 5 we prove Theorem 0.3 (stated as Theorem 5.1). 
Also, we use the incremental patching technique to show (see 5.3) that if $A_n \subset M_n$ are Cartan subalgebras in finite factors, with $\dim M_n \rightarrow \infty$,  
and $\Gamma_i$ are countable subgroups of the normalizer $\Cal N$ 
of $\text{\bf A}=\Pi_\omega A_n$ in $\text{\bf M}=\Pi_\omega M_n$,  
acting freely on  $\text{\bf A}$, with $H_i\subset \Gamma_i$ isomorphic amenable subgroups, then there exists $u\in \Cal N$ such that 
$uH_1u^*=H_2$ and such that the group generated by $u\Gamma_1u^*$ and $\Gamma_2$ is the amalgamated free product $\Gamma_1 *_H \Gamma_2$, 
where $H$ is the identification of $H_1, H_2$ via Ad$(u)$. Taking $M_n$ finite dimensional, this recovers  a result from [ES], [Pa], on the soficity of amalgamated free products 
of sofic groups over amenable subgroups and on the uniqueness of sofic embeddings of an amenable group. 

This paper was completed during my stay at the Jussieu Math Institute in Paris during the academic year 2012-2013. 
I want to gratefully acknowledge A. Connes, G. Pisier, G. Skandalis and S. Vassout, for their kind hospitality and support.

\heading 1. Preliminaries \endheading

\noindent
{\bf 1.1. Some generalities}. All von Neumann algebras $M$ considered in this paper are finite (in the sense of [MvN1]) 
and come equipped with a fixed faithful normal {\it trace} state, generically denoted $\tau$. 
We denote by $\Cal U(M)$  the group of unitary elements of $M$ and by $\Cal P(M)$ the set of projections of $M$. Recall that a von Neumann algebra is a 
{\it factor} if its center is reduced to the scalars. Recall that there exists a unique trace state on a finite factor ([D2]). 
A finite factor $M$ is either finite dimensional (in which case $M\simeq M_{n \times n}(\Bbb C)$ for some $n\geq 1$ 
with its unique trace state $\tau$ given by the normalized trace $tr=Tr/n$) 
or infinite dimensional. In this latter case, it is called a II$_1$ {\it factor}, and is characterized by the fact that 
the range of the trace on the set of projections satisfies $\tau(\Cal P(M))=[0,1]$. 

More generally, a finite von Neumann algebra splits as a direct sum $M=M_1 \oplus M_2$ 
with $M_1$ of {\it type} I (i.e. $M_1 \simeq \oplus_{n\geq 1} M_{n\times n}(\Bbb C) \otimes A_n$, where $A_n$ are abelian von Neumann algebras, possibly equal to 0)  
and $M_2$ of {\it type} II$_1$ (which by definition means $M_2$ has no type I summand).   

We denote by $\|x\|_2=\tau(x^*x)^{1/2}$, $x\in M$, the $L^2$ Hilbert-norm given by the trace. We denote by $L^2M$ 
the completion of $M$ in this norm. We  often view $M$ in its {\it standard representation}, acting on $L^2M$ by left multiplication.

We will also use the $L^1$ norm $\| \ \|_1$ on $M$, 
defined by  $\|x\|_1:=\tau(|x|)=\sup \{|\tau(xy)| \mid y\in M, \|y\|\leq 1\}$. We denote by 
$L^1M$ the completion of $M$ in the norm $\| \ \|_1$. Note that by the Cauchy-Schwartz inequality we have $\|x\|_1 \leq \|x\|_2$, 
while by the inequality $x^*x \leq \|x\| |x|$ we have $\|x\|^2_2 \leq \|x\|_1 \|x\|$.  

If $B \subset M$ is a von Neumann 
subalgebra, then $E_B: M \rightarrow B$ denotes the (unique) $\tau$-preserving conditional expectation of $M$ onto $B$, which is contractible 
in both the operatorial norm $\| \ \|$ and the above $L^p$-norms, $p=1,2$. If we view $M$ in its standard representation 
on $L^2M$, then the expectation $E_B$ is implemented by 
the orthogonal projection $e_B$ of $L^2M$ onto $L^2B\subset L^2M$ (viewed as the closure in the norm $\| \ \|_2$ of $B\subset M$), 
by $e_B x e_B = E_B(x)e_B$, $x\in M$.   

Given a von Neumann subalgebra $B \subset M$ and a set $X\subset M$, we say that $X$ is {\it perpendicular to} $B$ and write $X\perp B$ 
if $\tau(x^*b)=0$, $\forall x\in X$ and $b\in B$. 

A finite von Neumann algebra $(M, \tau)$ is {\it separable} if it is separable with respect to the Hilbert norm $\| \ \|_2$. 
Note that this condition is equivalent to the fact that $M$ is countably generated as a von Neumann algebra. More generally, if $X\subset M$ 
is a subspace, then $X$ is separable if it is separable with respect to the norm $\| \ \|_2$. 

The von Neumann algebra $M$ is {\it atomic} if $1_M=\Sigma_i e_i$ with $e_i\in M$ a family of mutually orthogonal 
minimal (or atomic) projections $e_i\in M$, i.e. with the property that $e_iMe_i=\Bbb Ce_i$.  $M$ is {\it diffuse} if it has no minimal (non zero) 
projection. Any abelian von Neumann algebra $A$ which is diffuse and separable is isomorphic to $L^\infty([0, 1])$ (or equivalently, to $L^\infty (\Bbb T)$). Moreover, if $A$ 
is endowed with a faithful normal state $\tau$, then the isomorphism $A \simeq L^\infty([0,1])$ can be taken so that to carry $\tau$ 
onto the integral $\int \cdot \ \text{\rm d}\mu$, where $\mu$ is the Lebesgue measure on $[0, 1]$. 

We will often consider maximal abelian $^*$-subalgebras (MASA) $A$ in  a finite von Neumann algebra $M$, i.e. abelian $^*$-subalgebras $A\subset M$ 
with $A'\cap M = A$. In such a case, 
we denote $\Cal N_{M}(A)=\{u\in \Cal U(M)\mid uAu^*=A\}$, the {\it normalizer} of $A$ in $M$. Following [FM], if the normalizer generates $M$ as a von 
Neumann algebra, we call $A$ a {\it Cartan subalgebra} in $M$. An isomorphism of Cartan inclusions 
$(A_0\subset M_0; \tau) \simeq (A_1\subset M_1; \tau)$ is a trace preserving isomorphism of $M_0$ onto $M_1$ carrying $A_0$ onto $A_1$. 

If $A_0\subset M_0$ is Cartan 
and $A_1\subset M_1$ is an arbitrary MASA, then 
a {\it Cartan embedding} (or simply an {\it embedding}) of $A_0\subset M_0$ into $A_1\subset M_1$ is a trace preserving embedding of $M_0$ into $M_1$ that carries 
$A_0$ into $A_1$ such that $M_0\cap A_1=A_0$, with the commuting square condition $E_{A_1}E_{M_0}=E_{A_0}$ satisfied (see 1.2 below for more on this condition), 
and such that  $\Cal N_{M_0}(A_0)\subset \Cal N_{M_1}(A_1)$. 

For various other general facts about finite von Neumann algebras, we refer the reader to the classic book [D2]. 

\vskip .05in 
\noindent
{\bf 1.2. Commuting squares of subalgebras}. 
Two von Neumann subalgebras $B_1, B_2\subset M$ are in {\it commuting square} position if the expectations $E_{B_1}, E_{B_2}$ commute (see Sec. 1.2 in [P2]). 
Note that if this is the case then we in fact have $E_{B_1}E_{B_2}=E_{B_2}E_{B_1}=E_{B_1\cap B_2}$. Also, 
for this to happen it is sufficient that $E_{B_1}(B_2)\subset B_1\cap B_2$. 

A typical example  when the commuting square condition is satisfied is the following: 
let $Q\subset P \subset M$ be von Neumann algebras; then $P$ and $Q'\cap M$ are in commuting square position  (see 1.2.2 in [P2]). 

We notice here an observation showing that in the statement of Theorem 0.1, we may equivalently take the space $X$ to be a separable 
von Neumann algebra making a commuting square with $\text{\bf Q}'\cap \text{\bf M}$, a fact that we will not use in the sequel but is good to keep in mind. 
See also (3.8 in [P12]) for a similar statement. 

\proclaim{Lemma} Let $N\subset M$ be a von Neumann subalgebra in the finite von Neumann algebra $M$. If $X\subset M$ is a separable subspace, 
then there exists a separable von Neumann subalgebra $P\subset M$ that contains $X$ and makes a commuting square with $N$.  
\endproclaim
\noindent
{\it Proof}. Let $P_0\subset M$  be the (separable) von Neumann algebra 
generated by $X$. We then denote by $B_1$ the von Neumann algebra generated by $E_N(P_0)$ and by $P_1$ the von Neumann algebra generated 
by $B_1$ and $P_0$. Note that $B_1\subset P_1$ are separable, with $B_1\subset N$ and $X\subset P_1$. More generally, we 
construct recursively an increasing sequence of inclusions of separable von Neumann algebras $B_n \subset P_n$, $n\geq 1$, 
by letting $B_n$ be the von Neumann algebra generated by $E_N(P_{n-1})$ and 
$P_n$ be the von Neumann algebra generated by $B_n$ and $P_{n-1}$. 

If we now define $B=\overline{\cup_n B_n}^w$ and $P=\overline{\cup_n P_n}^w$, then both algebras are separable and $B\subset P\cap N$, 
by construction. Moreover, we have $E_N(P_n)\subset B_{n+1}\subset P_{n+1}\subset P$, implying that $E_N(P)\subset B \subset P\cap N$. Thus, we actually have $E_N(P)=B=N \cap  P$, 
i.e.,  $N, P$ make a commuting square with $B=N\cap P$. 
\hfill 
$\square$

\vskip .05in 
\noindent
{\bf 1.3. Amenable algebras}. An important example of a (separable) II$_1$ factor is the {\it hyperfinite} II$_1$ {\it factor} $R$ of Murray and von Neumann ([MvN2]),  
defined as the infinite tensor product $(R, \tau)=\overline{\otimes}_k (M_{2 \times 2}(\Bbb C), tr)_k$. By [MvN2],  
$R$ is the unique {\it approximately finite dimensional} ({\it AFD}) separable II$_1$ factor (a separable finite von Neumann algebra 
algebra $(M,\tau)$ is AFD if there exists an increasing  sequence of finite dimensional von Neumann subalgebras $M_n\subset M$ such that $\cup_n M_n$ 
is dense in $M$ in the norm $\| \ \|_2$). 

By Connes' results in [C1], $R$ is in fact the unique {\it amenable} separable II$_1$ factor. 
Recall in this respect that a finite von Neumann algebra $(M, \tau)$ is called amenable if there exists a state $\varphi$ on $\Cal B(L^2M)$ 
that has $M$ (when viewed in its standard representation on $L^2M$) in its centralizer, $\varphi(xT)=\varphi(Tx)$, $\forall x\in M$, 
$\forall T \in \Cal B(L^2M)$, and such that $\varphi_{|M}=\tau$. Note that the latter condition is redundant in case $M$ 
is a factor, because $\varphi_{|M}$ is a trace and because of the uniqueness of the trace on factors. 
Connes Fundamental Theorem in [C1] 
actually shows that amenability is equivalent to the AFD property, for any finite von Neumann algebra. 

From all this, it follows that $R$ can be represented in many different ways, for instance as the {\it group measure space} 
II$_1$ factor $L^\infty(X)\rtimes \Gamma$, 
associated with a free  ergodic measure preserving action of a countable amenable group $\Gamma$ on a probability space $(X, \mu)$ ([MvN2]). 
When viewed this way, $R$ has $D=L^\infty(X)$ as a natural Cartan subalgebra. By [CFW], [OW] the Cartan subalgebra of $R$ is in fact unique, up 
to conjugacy by an automorphism of $R$. We may thus represent $D\subset R$ as the infinite tensor product  $\overline{\otimes}_k (D_2)_k\subset 
\overline{\otimes} (M_{2 \times 2}(\Bbb C))_k$, where 
$D_2$ is the diagonal subalgebra in $M_{2 \times 2}(\Bbb C)$. 

More generally, by [CFW], if $A_0\subset R_0$ is a Cartan subalgebra in an amenable separable finite von Neumann algebra $R_0$, 
then there exists an increasing sequence of finite dimensional Cartan inclusions $(A_{0,n}\subset R_{0,n}) \subset (A_0 \subset R_0)$  
(with Cartan embeddings, as defined before) such that $\overline{\cup_n A_{0,n}}^w=A_0\subset R_0 = \overline{\cup_n R_{0,n}}^w$.

\vskip .05in
\noindent
{\bf 1.4. Local quantization relative to subalgebras}. We recall here a result from [P1], [P5], showing  that 
if $Q\subset M$ are II$_1$ von Neumann algebras, then one can ``simulate'' the expectation 
onto the commutant $Q'\cap M$ by ``squeezing'' with appropriate projections in $Q$, a phenomenon called ``local quantization'' 
in [P5]:    

\proclaim{Theorem} $1^\circ$ Let $M$ be a finite von Neumann algebra and  $Q\subset M$ a von Neumann subalgebra. Given any  
finite set $F\subset M\ominus Q\vee (Q'\cap M)$ and any $\varepsilon > 0$, there exists a projection $q\in Q$ 
such that $\|q x q\|_1 < \varepsilon \tau(q)$, $\forall x\in F$. 

$2^\circ$ Let $Q\subset M$ be an inclusion of $\text{\rm II}_1$ von Neumann algebras. 
Given any finite set $X\subset M$ and any $\varepsilon > 0$, there exists a projection $q\in Q$ 
such that $\|qxq-E_{Q'\cap M}(x)q\|_1 < \varepsilon \tau(q)$, $\forall x\in X$. Moreover, $q$ can be taken so that to 
have scalar central trace in $Q$. 
\endproclaim  

\noindent 
{\it Proof}. Part $1^\circ$ is already proved in [P1] (see also Theorem 3.6 in [P12]),  while part $2^\circ$ is (Theorem A.1.4  in [P5]). 

\hfill 
$\square$ 

\vskip .05in  

\noindent
{\bf 1.5. A criterion for non-conjugacy of subalgebras}. Let $Q, P \subset M$ be von Neumann subalgebras of the finite von Neumann 
algebra $M$. Following [P10], we say that {\it a corner of $Q$ can be embedded into $P$ 
 inside $M$} and write $Q\prec_M P$ if the following condition holds true:  there exist non-zero projections $p\in P$, $q\in Q$, a unital
isomorphism $\psi: qQq \rightarrow pPp$ (not necessarily onto)
and a partial isometry $v\in
M$ such that $vv^*\in (qQq)'\cap qMq$, $v^*v \in \psi(qQq)'\cap pMp$,  
$xv = v\psi(x), \forall x\in qQq$, and $x \in qQq$, $xvv^*=0$, implies $x=0$. 

In this paper we will actually consider cases  when the above condition is 
not satisfied. We recall from (2.1 in [P10]) a useful necessary and sufficient criterion for this to happen:

\proclaim{Theorem} Let $M$ be a finite
von Neumann algebra and $P,Q\subset M$ von Neumann subalgebras. For each $q\in \Cal P(Q)$, 
fix $\Cal U_q\subset \Cal U(qQq)$ a subgroup generating $qQq$ as a von Neumann algebra. 
Then $Q\not\prec_M P$ if and only if the following condition holds true:  
\vskip .05in 
\noindent
$(1.5.1)$ Given any $q\in \Cal P(Q)$ and any separable subspace $X\subset M$
there exists a sequence of unitary elements $u_n\in \Cal U_q$ such that $\lim_n \|E_P(xu_ny)\|_2$ $ =0$, $\forall x, y  \in X$. 
\vskip .05in 
\endproclaim

\vskip .05in
\noindent
{\bf 1.6. Ultraproducts of algebras}. We fix once for all an (arbitrary) 
free ultrafilter $\omega$ on $\Bbb N$.  If $M_n$, $n\geq 1$, is a sequence of finite von Neumann algebras then, we denote by 
$\Pi_\omega M_n$ their $\omega$-{\it ultraproduct}, i.e., the finite von Neumann algebra obtained as the quotient of $\oplus_n M_n$ by its ideal 
$\Cal I_\omega=\{ (x_n) \mid \lim_\omega \tau(x_n^*x_n) = 0\}$, endowed with the trace $\tau(y) =\lim_\omega \tau(y_n)$, where $(y_n)_n \in \oplus_n M_n$ 
is in the class $y \in \oplus_n M_n/\Cal I_\omega$ ([Wr]). Recall that if $M_n$ are factors and $\dim M_n \rightarrow \infty$, 
then $\Pi_\omega M_n$ is a II$_1$ factor ([Wr]) and it is non-separable ([F]). 

If $Q_n\subset M_n$ are von Neumann subalgebras, $n\geq 1$, then the ultraproduct 
$\Pi_\omega Q_n$ identifies naturally to a von Neumann subalgebra in $\Pi_\omega M_n$ and its {\it centralizer} (or commutant) in $\Pi_\omega M_n$ is given by the formula   
$(\Pi_\omega Q_n)' \cap \Pi_\omega M_n=\Pi_\omega (Q_n'\cap M_n)$ 
(see e.g. [P1]). 

If $M$ is a finite von Neumann algebra, then $ M^\omega$ denotes its $\omega$-{\it ultrapower}, 
i.e. the ultraproduct of infinitely many copies of $M$. Note that $M$ naturally embeds into $M^\omega$, as 
the von Neumann subalgebra of constant sequences, and that if $M$ is a II$_1$ factor then $M^\omega$ is a (non-separable by [F]) II$_1$ factor. 

\vskip .05in
\noindent
{\bf 1.7. Centralizers of countable sets in ultraproducts}.  Let $S=\{b_n\}_n$ be a countable subset in the ultrapower $R^\omega$  of the hyperfinite 
II$_1$ factor $R$ and let $b_n=(b_{n,m})_m$ be representations of each of its elements with $b_{n,m}\in R=\overline{\otimes}_k (M_{2 \times 2}(\Bbb C))_k 
=\overline{ \cup_n M_n}^w$, where $M_n$ is the tensor product of the first $n$ copies of $M_{2\times 2}(\Bbb C)$. Thus, we may assume that for each $m$, 
$\{b_{n,m}\}_{n\leq m}\subset M_{k_m}$, for a large enough $k_m$. Then we have $b_n \in \Pi_\omega  M_{k_m} \subset R^\omega$, $\forall n$, 
viewed as a subalgebra of $R^\omega$.  But then the ultraproduct subalgebra $\Pi_\omega (M_{k_m}'\cap R)\simeq R^\omega$ 
commutes with the set $\{b_n\}_n$. This shows that the centralizer of any separable von Neumann subalgebra $B$ of $R^\omega$ is a type II$_1$ 
von Neumann algebra without separable direct summands. 

More generally, the same argument shows that if $\text{\bf M}=\Pi_\omega M_n$ is an ultraproduct of arbitrary McDuff II$_1$ factors $M_n$ (i.e., for which we have 
$M_n\simeq M_n \overline{\otimes} R$, see [McD]), then 
the centralizer of any separable subalgebra $B \subset \text{\bf M}$ is of type II$_1$ with no separable direct summands.

However, for general ultraproducts $\Pi_\omega M_n$ and ultrapowers $M^\omega$, we may have countable (or even finite) subsets $S$ that have trivial 
centralizer: for instance, if $M$ is a separable non-Gamma II$_1$ factor ([MvN2]), then $M$ is countably generated and $M'\cap M^\omega=\Bbb C1$ (by [McD]). 
This is the case if $M=L(\Bbb F_n)$, with $\Bbb F_n$ the free group with $2\leq n\leq \infty$ generators (cf. [MvN2]), or if 
$M=L(\Gamma)$ with $\Gamma$ an ICC group with the property (T) of Kazhdan (for example, $\Gamma=PSL(n, \Bbb Z)$, $n \geq 3$). 
Similarly, by results in [Be], it follows that if for some fixed $n\geq 3$ 
we take $(\pi_m, \Cal H_m)$ to be any sequence of finite dimensional irreducible representations of $\Gamma=PSL(n,\Bbb Z)$ so that 
$k_m=\dim \Cal H_m \rightarrow \infty$, then the von Neumann subalgebra $M$ generated by $\{(\pi_m(g))_m  \mid g\in \Gamma\}$ 
in the ultraproduct II$_1$ factor $\Pi_\omega M_{k_m\times k_m}(\Bbb C)$ is isomorphic to the group factor $L(\Gamma)$ and 
has trivial relative commutant. 

The following result shows that in fact the centralizer of any separable von Neumann subalgebra 
$P$ of an arbitrary ultraproduct II$_1$ factor 
$\text{\bf M}:=\Pi_\omega M_n$, coming from a sequence of finite factors $M_n$ with 
$\dim M_n \rightarrow \infty$, splits as the direct sum of an atomic von Neumann algebra and a diffuse von Neumann algebra with only non-separable direct 
summands.

\proclaim{Theorem} If $P$ is a separable von Neumann subalgebra of $\text{\bf M}$ then 
$P'\cap \text{\bf M}=B_0\oplus B_1$, with $B_0$ atomic and $B_1$ diffuse and having no  
separable direct summand $($even more: any MASA of $B_1$ has only non-separable direct summands$)$.  
\endproclaim
\noindent
{\it Proof}. Denote $Q=P'\cap \text{\bf M}$ and let $z\in \Cal Z(Q)$ be the maximal central projection with 
the property that $Qz$ is diffuse. We have to prove that $Qz'$ is non-separable for any central projection $z'\in \Cal Z(Q)z$. 
By replacing $P\subset \text{\bf M}$ by $Pz\subset z\text{\bf M}z$, we may clearly assume $z=1$. 

Assuming by contradiction that $Q$ has separable direct summands, we may further reduce with 
the maximal central projection $z_0$ in $Q$ with the property that $Qz_0$ is separable to actually assume, 
by contradiction, that $P\subset \text{\bf M}$ is separable with $Q=P'\cap \text{\bf M}$ diffuse and separable. 

Let $\{b_n\}_n\subset P$ be a countable subset of the unit ball of $P$, dense in the Hilbert norm $\| \ \|_2$. Let $b_n=(b_{n,m})_m$ 
be representations of $b_n$ with $b_{n,m}\in M_m$, $\|b_{n,m}\|\leq \|b_n\|$, $\forall n,m$. 
Let also $u \in Q$ be a Haar unitary generating a maximal abelian $^*$-subalgebra $A_0$ of $Q$, 
and let $u=(u_m)_m$ be a representation of $u$ with $u_m\in \Cal U(M_m)$, $\forall m$.  

The fact that $u$ belongs to $Q=\{b_n\}_n'\cap \text{\bf M}$ translates into the condition
$$
\underset m \rightarrow \omega \to \lim \|[b_{k,m}, u_m]\|_2 = 0, \forall k \geq 1,  \tag 1
$$
while the fact that $u$ is a Haar unitary amounts to the condition
$$
\underset m \rightarrow \omega \to \lim \tau(u_m^j) = 0, \forall j\neq 0. \tag 2
$$

Let $V_n$ denote the set of all $m\in \Bbb N$ with the property that 
$$
\|[b_{k,m}, u_m]\|_2 < 2^{-n}, |\tau(u_m^j)|<2^{-n}, \forall 1\leq k \leq n, 1\leq |j| \leq 2n. \tag 3
$$

If we identify $\ell^\infty\Bbb N$ with the algebra $C(\Omega)$ of continuous functions on its spectrum $\Omega$ (via the GNS representation), 
and we view $\omega$ as a point in $\Omega$, 
then by $(1)$ and $(2)$ it follows that $V_n$ correspond to an open-closed neighborhoods of $\omega \in \Omega$.  
Let now $W_n$, $n\geq 0$,  be defined 
recursively as follows: $W_0=\Bbb N$ and $W_{n+1}=W_n \cap V_{n+1} \cap \{n \in \Bbb N \mid n > \min  W_n\}$. Note that, 
with the same identification as before, $W_n$ correspond to  
a strictly decreasing sequence of neighborhoods of $\omega$. 

Noticing that the sets $\{W_{n-1}\setminus W_n\}_{n\geq 1}$ form a partition of $\Bbb N$, 
we define $v=(v_m)_m$ by letting $v_m=u_m^n$ for $m\in W_{n-1} \setminus W_n$. 
Since $v_m\in \Cal U(M_m)$, it follows that $v$ is a unitary element in $\text{\bf M}$. By the first relation 
in $(3)$, if $m\in W_{n-1}\setminus W_n$ then 
$$
\|[b_{k,m},v_m]\|_2 = \|[b_{k,m},u_m^n]\|_2 \leq \Sigma_{j=0}^{n-1} \|u_m^j[b_{k,m},u_m]u_m^{n-j-1}\|_2 \leq n2^{-n}, \tag 4
$$
for all $1\leq k \leq n$, while by the second relation in $(3)$ we have 
$$
|\tau (v_m u_m^j)| < 2^{-n} \tag 5 
$$
for all $1\leq |j| \leq n$. 

But then $(4)$ implies $v\in \{b_n\}_n'\cap \text{\bf M}=P'\cap \text{\bf M}=Q$, while by $(5)$ we have $\tau(vu^j)=0$, 
for all $j\neq 0$, i.e. $v\in Q$ is perpendicular to the maximal abelian $^*$-subalgebra $A_0=\{u\}''$ of $Q$ generated by $u\in Q$.  
Since by construction we have $uv=vu$, this shows that at the same time we have $v\in \{u\}'\cap Q=A_0$ 
and $v\perp A_0$,  a contradiction. This also shows the stronger form of the statement.
\hfill $\square$

\heading 2. Bicentralizer characterizations of amenability \endheading

\proclaim{2.1. Theorem} $1^\circ$ Let $M_n$ be a sequence of finite factors with $\text{\rm dim} M_n \rightarrow \infty$ and  
denote $\text{\bf M}=\Pi_\omega M_n$. If $B\subset \text{\bf M}$ is a separable amenable von Neumann subalgebra, then 
$(B'\cap \text{\bf M})'\cap  \text{\bf M}=B$. Moreover, $B'\cap \text{\bf M}$ is 
of type $\text{\rm II}_1$ and has only non-separable direct summands. 

$2^\circ$ If $R$ denotes the hyperfinite  $\text{\rm II}_1$ factor  then $(R'\cap R^\omega)'\cap R^\omega=R$.  
\endproclaim 
\noindent
{\it Proof}. Part $2^\circ$ is just a particular case of part $1^\circ$, so we only need to prove $1^\circ$. 
By Connes' Theorem ([C1]), since $B$ is amenable and separable, it is approximately finite dimensional, so 
$B = \overline{ \cup_n B_n}^w$, for some increasing sequence of finite dimensional 
von Neumann subalgebras $B_n\subset B$. Note that $B'\cap \text{\bf M}=\cap_n (B_n'\cap \text{\bf M})$ and that for each $n$ 
we have $(B_n'\cap \text{\bf M})'\cap \text{\bf M}=B_n$ (in fact, it is trivial to see that given any inclusion of von Neumann algebras 
$\Cal N \subset \Cal M$ with $\dim \Cal N < \infty$ and $\Cal M$ a factor, we have $(\Cal N'\cap \Cal M)'\cap \Cal M=\Cal N$). We 
first need to prove the following: 

\vskip .05in 

{\it Fact}. Let $P \subset M$ be an inclusion of finite von Neumann algebras. Let $x\in M\ominus (P'\cap M)$, $x\neq 0$, and $\varepsilon > 0$. 
There exists a unitary element $u\in P$ such that $\Re \tau(x^*uxu^*) < \varepsilon \|x\|_2^2$. 

\vskip .05in

To prove this, let $K_x$ denote the weak closure of the convex set $\text{\rm co} \{uxu^* \mid u\in \Cal U(P) \}$ and  
note right away that $\|y\|\leq \|x\|$ and $\|y\|_2 \leq \|x\|_2$, $\forall y \in K_x$. Thus, $K_x$ is a weakly closed bounded 
subspace in both $M$ and $L^2M$. In particular, there exists a unique element $y_0 \in K_x$ 
of minimal Hilbert-norm: $\|y_0\|_2 = \min \{ \|y\|_2 \mid y\in K_x\}$. Since $K_x$ is Ad$(\Cal U(P))$-invariant  
(because it is the weak closure of the Ad$(\Cal U(P))$-invariant set $\text{\rm co} \{uxu^* \mid u\in \Cal U(P)\}$) and since 
$\|uy_0u^*\|_2=\|y_0\|_2$, by the uniqueness of $y_0$ it follows that $uy_0u^*=y_0$, $\forall u \in \Cal U(P)$. Thus, 
$uy_0=y_0u$, $\forall u\in \Cal U(P)$. By taking linear combinations of $u$, this implies $y_0 \in P'\cap M$. 
But by its construction, the entire $K_x$ lies in $M\ominus (P'\cap M)$. Thus, $y_0$ is both in $P'\cap M$ and perpendicular to it,  
implying that $y_0=0$, i.e. $0\in K_x$. 

Assuming now that we have $\Re \tau(x^*uxu^*) \geq \varepsilon \|x\|_2^2$, for all $u \in \Cal U(P)$, 
by taking convex combinations over $u\in \Cal U(P)$  and then weak closure, it follows that $\Re \tau(x^*y) \geq \varepsilon \|x\|_2^2$, for all $y \in P$. 
In particular, $0=\Re \tau(x^*y_0) \geq \varepsilon \|x\|_2^2$, forcing $x=0$, a contradiction. This ends the proof of the above {\it Fact}. 

\vskip .05in

Denote for simplicity $Q=B'\cap \text{\bf M}$ and note that $B\subset Q'\cap \text{\bf M}$. Assume there exists 
$x\in Q'\cap \text{\bf M}$ with $x\perp B$.  In particular $x \perp B_n = (B_n'\cap \text{\bf M})'\cap \text{\bf M}$. 
By applying the {\it Fact} to the inclusion $B_n'\cap \text{\bf M} \subset \text{\bf M}$ and the element $x$, it follows that there exists a unitary 
element $u_n \in B_n'\cap \text{\bf M}$ such that $\Re \tau(x^*u_nxu_n^*) < 2^{-n}$, $\forall n$. 

Let $\{e^n_k\}_k\subset B_n$ denote the (finite) pseudogroup of all partial isometries in $B_n$  that can be obtained as a sum of 
elements from a given matrix unit of $B_n$, and which we take so that  $\{e^n_i\}_i$ is a subset of $\{e^{n+1}_j\}_j$, $\forall n$.  
Let $e^n_k=(e^n_{k,m})_m$, with $e^n_{k,m}\in M_m$ chosen so that $\|e^n_{k,m}\|\leq \|e^n_k\|$ 
and $\{e^n_{i,m}\}_i \subset  \{e^{n+1}_{j,m}\}_j$ for all $n , m$. Let also $u_n=(u_{n,m})_m$, with $u_{n,m}\in \Cal U(M_m)$. Then the above properties 
translate into 
$$
\underset m \rightarrow \omega \to \lim \|[u_{n,m}, e^n_{k,m}]\|_2 = 0, \underset m \rightarrow \omega \to \lim \Re \tau(x_m^*u_{n,m}x_mu_{n,m}^*) < 2^{-n}, \tag 1 
$$
for all $k$ and all $n$, where $x=(x_m)_m$ with $x_m \in M_m$. 

Let $V_n$ denote the set of all $m\in \Bbb N$ with the property that 
$$
\|[u_{n,m}, e^n_{k,m}]\|_2 < 2^{-n}, \Re \tau(x_m^*u_{n,m}x_mu_{n,m}^*) < \|x\|_2^2/2, \forall k. \tag 2 
$$

By $(1)$,  it follows that $V_n$ corresponds to an open-closed neighborhood of $\omega$ in the spectrum $\Omega$ 
of $\ell^\infty \Bbb N$, under the identification 
$\ell^\infty \Bbb N = C(\Omega)$. Let now $W_n$, $n\geq 0$,  be defined 
recursively as follows: $W_0=\Bbb N$ and $W_{n+1}=W_n \cap V_{n+1} \cap \{n \in \Bbb N \mid n > \min  W_n\}$. Note that, 
with the same identification as before, $W_n$ correspond to  
a strictly decreasing sequence of neighborhoods of $\omega$. 
Define $v=(v_m)_m$ by letting $v_m=u_{n,m}$ for $m\in W_{n-1} \setminus W_n$. 
Since $v_m\in \Cal U(M_m)$, it follows that $v$ is a unitary element in $\text{\bf M}$, while by the first relation 
in $(2)$ and the fact that $\{e^n_{i,m}\}_i \subset \{e^{n+1}_{j,m}\}_j$ it follows that $v\in \cap_n B_n'\cap \text{\bf M}=B'\cap \text{\bf M}=Q$. 
By the second relation in $(2)$, we also have $\Re \tau(x^*vxv^*) \leq \|x\|_2^2/2$. But  $x \in Q'\cap \text{\bf M}$ by our assumption, 
thus $vxv^*=x$, giving $\tau(x^*vxv^*)=\|x\|_2^2$, a contradiction. 

If $Q=Qz + Q(1-z)$ with $z$ a non-zero central projection of $Q$ and $Qz$ separable, then by the bi-commutant property 
we have $z\in B$ and by Theorem 1.7,  $Qz$ is atomic. Thus,  $Bz=(Qz)'\cap z\text{\bf M}z$ would follow non-separable, a 
contradiction. 

Assume now that $Q=Qz+Q(1-z)$ with $z\in \Cal P(\Cal Z(Q))$ such that $Qz$ is type I. By the bi-commutation relation, it follows  
again that $z\in B$ and that $Bz=(Qz)'\cap z\text{\bf M}z$ is non-separable (because the commutant of any abelian 
von Neumann subalgebra of $\text{\bf M}$ is non-separable, by 4.3 in [P1], or 2.3 in [P12]). 
\hfill $\square$

\proclaim{2.2. Theorem} Let $A_n \subset M_n$ be a sequence of MASAs in finite factors and   
denote $\text{\bf A}=\Pi_\omega A_n \subset \Pi_\omega M_n=\text{\bf M}$, $\Cal N=\Cal N_{\text{\bf M}}(\text{\bf A})$. 

$1^\circ$ If $H\subset \Cal N$ is a countable amenable 
subgroup, then $(H'\cap \text{\bf A})'\cap \text{\bf M}=\text{\bf A}\vee H$. 

$2^\circ$ Assume the MASAs $A_n \subset M_n$ are Cartan. Let $R_0\subset \text{\bf M}$ be a separable amenable 
von Neumann subalgebra such that $D_0=R_0 \cap \text{\bf A}$ is a Cartan subalgebra in $R_0$ and such that   
$(D_0\subset R_0)$ is Cartan embedded into $(\text{\bf A}\subset \text{\bf M})$, in the sense of $1.1$.  
Then $(\Cal N_{R_0}(D_0)'\cap \Cal N)'\cap \Cal N=  \Cal N_{R_0}(D_0)$. 
Moreover, if $D_1\subset R_1$ is another Cartan inclusion 
which is Cartan embedded 
into $\text{\bf A}\subset \text{\bf M}$, then given any isomorphism $\rho: (D_0\subset R_0; \tau)\rightarrow (D_1\subset R_1; \tau)$, 
there exists $u\in \Cal N$ such that $\text{\rm Ad}(u)=\rho$ on $R_0$.

$3^\circ$ With the same assumptions and notations as in $2^\circ$ above, let $\text{\bf A}_0=R_0'\cap \text{\bf A}$ and $\Cal N_0=\Cal N_{R_0}(D_0)'\cap \Cal N$. Then $\text{\bf A}_0$ is maximal abelian in $R_0'\cap \text{\bf M}$, $\Cal N_0$ coincides with the normalizer of $\text{\bf A}_0$ in $R_0'\cap \text{\bf M}$ and $\text{\bf M}_0 = \text{\bf A}_0 \vee \Cal N_0$  satisfies 
$\text{\bf M}_0'\cap \text{\bf M}=R_0$. 
\endproclaim 
\noindent
{\it Proof}. $1^\circ$ Let first  $\{e^n_j\}_j$ be an increasing sequence of finite partitions in $\Cal P(\text{\bf A})$ 
such that $\lim_n \|\Sigma_j e^n_j u e^n_j - E_{\text{\bf A}}(u)\|_2 = 0$, $\forall u\in H$ (e.g., by [P1], or 3.6 in [P12]). If we denote by $A_0$ the von Neumann 
subalgebra of $\text{\bf A}$ generated by $\cup_{u\in H}u\{e^n_j\mid j, n\}u^*$ and $R_0 = A_0 \vee H$, then $H$ normalizes $A_0$,  
$A_0$ is a Cartan subalgebra of $R_0$ and $\text{\bf A}\vee H = \text{\bf A}\vee R_0$. In particular, $H'\cap \text{\bf A}=R_0'\cap \text{\bf A}$. 
Moreover, since $H$ is amenable, $R_0$ 
follows amenable so by ([CFW], [OW]) there exists an increasing sequence of finite pseudogroups of partial isometries $\Cal G_n=\{e^n_j\}_j$, normalizing $A_0$  
(and $\text{\bf A}$ as well), with source and targets either equal or mutually orthogonal, for each $n \geq 1$, and such that $\{e^n_j \mid j, n\}$ 
generate $R_0$. 

It is then trivial to see that $H'\cap \text{\bf A} = \cap_n (\Cal G_n'\cap \text{\bf A})$ and 
$(\Cal G_n'\cap \text{\bf A})'\cap \text{\bf A}=\Cal G_n \vee \text{\bf A}$, $\forall n$. Then the rest of the proof 
proceeds with a ``diagonalization'' argument, 
exactly as at the end of the proof of Theorem 2.1.  

$2^\circ$ The proof of this part is similar to the one of $2.1.1^\circ$ and of $2.2.1^\circ$ above. Indeed, the statement obviously 
holds true for $D_0\subset R_0$ finite dimensional. Then for general $D_0\subset R_0$  one takes $\Cal G_n$ as in the proof of $2.2.1^\circ$ and 
one denotes  
by $D(n)\subset R(n)$ the associated (finite dimensional) Cartan inclusion. Noticing that $\Cal N_{R_0}(D_0)'\cap \Cal N = \cap_n \Cal N_{R(n)}(D(n))'\cap \Cal N$,   
one then combines the finite dimensional case  with a diagonalization argument, as in the proof of 2.1. 

$3^\circ$ Note first that $\Cal N_0$ normalizes $\text{\bf A}_0$. Indeed, if $a_0\in \text{\bf A}$ commutes with $R_0$ and $u\in \Cal N_0$, then $ua_0u^* \in \text{\bf A}$ and it commutes 
with $R_0$ (because both $a_0$ and $u$ commute with $R_0$). 

To see that $\text{\bf A}_0$ is a MASA in $R_0'\cap \text{\bf M}$, note that by part $2^\circ$ above we have 
$\text{\bf A}_0'\cap \text{\bf M}=\text{\bf A}\vee \Cal N_{R_0}(D_0)$. Thus 
$$
\text{\bf A}_0'\cap (R_0'\cap \text{\bf M})=R_0'\cap (\text{\bf A}\vee R_0)=R_0'\cap (D_0'\cap 
 (\text{\bf A}\vee R_0))=R_0'\cap \text{\bf A}=\text{\bf A}_0, 
$$
where we have used the fact that $E_{D_0'\cap \text{\bf M}}(R_0)=D_0\subset \text{\bf A}$. This also shows that $(\text{\bf A}_0\vee D_0)'\cap \text{\bf M}=\text{\bf A}$. 
If now $u\in R_0'\cap \text{\bf M}$ is a unitary that normalizes $\text{\bf A}_0$, then $u$ commutes with $D_0$ so it normalizes $\text{\bf A}_0\vee D_0$, 
and thus also its commutant $\text{\bf A}$, i.e., $u\in \Cal N \cap R_0'=\Cal N_0$. 

Finally, by part 2$^\circ$ above, we have  
$\text{\bf M}_0'\cap \text{\bf M}=(\text{\bf A} \vee R_0)\cap \Cal N_0' = R_0$. 
\hfill 
$\square$

\vskip .05in 
\noindent
{\bf 2.3. Some remarks and open problems.} $1^\circ$ It is well known (and trivial to show) 
that if $M_n$ is a sequence of finite factors with $\dim M_n \rightarrow \infty$ 
and $(B, \tau)$ is a finite separable AFD von Neumann algebra, then there exists a trace preserving embedding 
$\theta_0: B \hookrightarrow \text{\bf M}:=\Pi_\omega M_n$ and that given any other such trace preserving embedding $\theta_1:B \hookrightarrow \text{\bf M} $, 
there exists a 
unitary element $u\in \text{\bf M}$ such that $\theta_1(b)=u\theta_0(b)u^*$, $\forall b\in B$. In particular, any two copies 
of $(B, \tau)$ in $\text{\bf M}$ are unitary conjugate. By Connes' theorem [C1], this means that the same  holds true for any 
finite, separable, amenable $B$. 

Moreover, by a result of K. Jung in [J], the converse is also true: if a finite separable 
von Neumann algebra $(B, \tau)$ has a unique (up to unitary conjugacy) embedding  into either an  
ultraproduct $\Pi_\omega M_{n \times n}(\Bbb C)$ or in $R^\omega$, then $B$ is amenable 
(see [J]). In fact, by a result of N. Brown in [B] (see also Ozawa's Appendix 8.1 in that paper), if $B\subset R^\omega$ is non-amenable, then there exist 
uncountably many non-conjugate copies of $B$ in $R^\omega$. 

Since given any ultraproduct II$_1$ factors $\text{\bf M}=\Pi_\omega M_n$, all embeddings $B \hookrightarrow \text{\bf M}$ 
of a given separable amenable finite von Neumann algebra are unitary conjugate in $\text{\bf M}$, 
it seems interesting to investigate the converse in this general setting: is it true that if $B\subset \text{\bf M}$ is a separable non-amenable 
von Neumann algebra of an arbitrary utraproduct II$_1$ factor, 
then there exist ``many'' non-conjugate copies of $B$ in $\text{\bf M}$?  (I am grateful to N. Ozawa for pointing out to me 
that the answer to this problem is not known; see [FHS] for related considerations.)

On the other hand, related to Theorem 2.1 above, we propose the following new characterization of amenability for 
separable finite von Neumann algebras:

\vskip .05in 
\noindent 
{\bf (2.3.1)}  {\it Conjecture}: Let $P$ be a separable von Neumann subalgebra 
of an ultraproduct II$_1$ factor $\text{\bf M}$ (notably, of $\text{\bf M}=R^\omega$, or of  
$\text{\bf M}=\Pi_\omega M_{n \times n}(\Bbb C)$). If 
the bicentralizer condition $(P'\cap \text{\bf M})'\cap \text{\bf M}=P$ is satisfied, then $P$ is amenable. 
In particular, if $M$ is a separable II$_1$ factor such that $(M'\cap M^\omega)'\cap M^\omega=M$, then $M \simeq R$. 

\vskip .05in 

Note that for a separable von Neumann subalgebra $P$ of an ultraproduct II$_1$ factor $\text{\bf M}$, being equal to its bicentralizer 
is equivalent to being equal to the centralizer of some $^*$-subalgebra of $\text{\bf M}$. Thus, conjecture $(2.3.1)$ is equivalent to the following: 

\vskip .05in
\noindent
{\bf (2.3.1')}  {\it Conjecture}: Let $P$ be a separable von Neumann subalgebra 
of an ultraproduct II$_1$ factor $\text{\bf M}$. If $P$ 
is the centralizer of a von Neumann subalgebra $\text{\bf Q}\subset \text{\bf M}$, i.e., $P= \text{\bf Q}'\cap \text{\bf M}$, 
then $P$ is necessarily amenable.  

\vskip .05in 

Indeed, one clearly has that $(2.3.1')$ implies $(2.3.1)$. Assume in turn that $(2.3.1)$ holds true. Let $\text{\bf Q}\subset \text{\bf M}$ 
be so that  $P=\text{\bf Q}'\cap \text{\bf M}$ is separable and denote $\tilde{\text{\bf Q}}=P'\cap \text{\bf M}$. Then we still have  
$\tilde{\text{\bf Q}}'\cap \text{\bf M}=P$, so $P$ satisfies the bicentralizer condition and it is separable, thus $P$ is amenable.

Note also that the bicentraliser condition  $(M'\cap M^\omega)'\cap M^\omega = M$ for a separable II$_1$ factor $M$, implies that $M$ must be McDuff ([McD]), 
i.e., it splits off the hyperfinite II$_1$ factor (or else $M'\cap M^\omega$ is abelian, implying that the bicentralizer is non-separable), but that it cannot be of the form $N\overline{\otimes}R$, with $N$ non-Gamma ([MvN2]). 
More generally, if $M$ has a II$_1$ von Neumann subalgebra $N\subset M$ satisfying the spectral gap condition $N'\cap M^\omega = (N'\cap M)^\omega$ ([P11]),  then $M$ 
cannot satisfy the bicentralizer condition $(M'\cap M^\omega)'\cap M^\omega=M$. Indeed, this is because 
taking bicentralizer is an  operation preserving  inclusions of algebras, and thus 
the bicentralizer of $M$ in $M^\omega$ contains the bicentralizer of $N$ in $M^\omega$, which is equal to $((N'\cap M)^\omega)'\cap M^\omega = N^\omega$. But the latter  
is non-separable, so it cannot be contained in $M$, which is separable.  Finally, note  that if $M$ is McDuff, then given any separable $^*$-subalgebra $B\subset M^\omega$, 
its centralizer $B'\cap M^\omega$ is of type II$_1$. More precisely, if $R=\overline{\otimes}_n (M_{2 \times 2}(\Bbb C))_n$, then there exists a sufficiently fast growing $k_n \rightarrow \infty $ 
such that if we denote $R_n=\overline{\otimes}_{m \geq k_n} (M_{2 \times 2}(\Bbb C))_m$, then $B'\cap M^\omega$ contains $\Pi_\omega R_n$.

\vskip .05in

$2^\circ$ Since by ([CFW]), any Cartan inclusion $A_0\subset M_0$ with $M_0$ separable amenable finite von Neumann algebra is a limit of an increasing 
sequence of finite dimensional Cartan inclusions (see 1.3), it follows that any 
isomorphism between two embeddings of $A_0\subset M_0$ into an ultraproduct inclusion 
$\text{\bf A}\subset \text{\bf M}$ is  implemented by a unitary element in $\Cal N_{\text{\bf M}}(\text{\bf A})$. Indeed, this is clear for finite 
dimensional $A_0\subset M_0$, and the general case follows by a diagonalization procedure.  

If in turn $A_0\subset M_0$ is a Cartan subalgebra with $M_0$ non-amenable, and $A_0\subset M_0$ is embeddable into 
an ultraproduct $\text{\bf A} \subset \text{\bf M}$ which is  
either of the form $\Pi_\omega D_n \subset \Pi_\omega M_{n\times n}(\Bbb C)$, or of the form $D^\omega \subset R^\omega$, then 
any two copies of $A_0\subset M_0$ into $\text{\bf A}\subset \text{\bf M}$ that are conjugate by a unitary in $\Cal N_{\text{\bf M}}(\text{\bf A})$  
will in particular have the corresponding copies of $M_0$ unitary conjugate in $\text{\bf M}$. The procedure of constructing 
``many'' non-conjugate embeddings of a non-amenable $M_0\subset \text{\bf M}$ starting from an initial embedding of $M_0$ in the proof of (8.1 of  [B]), is easily seen to 
actually give embeddings of $A_0\subset M_0$ into $\text{\bf A}\subset \text{\bf M}$, once the initial embedding of $M_0$ 
is in fact a Cartan embedding of $A_0\subset M_0$ into $\text{\bf A}\subset \text{\bf M}$. Thus, (8.1 in [B]) also implies that there exist 
uncountably many non-conjugate embeddings of $A_0\subset M_0$ into $\text{\bf A}\subset \text{\bf M}$. 
Altogether, this gives an analogue for Cartan inclusions (equivalently, for countable equivalence relations [FM]), of K. Jung's characterization 
of amenability in [J],  by a ``unique embedding'' - type property. 

Part 2$^\circ$ of Theorem 2.2 above suggests that, for a separable Cartan inclusion $A_0\subset M_0$ embedded into 
an ultraproduct of Cartan inclusions $\text{\bf A}\subset \text{\bf M}$, 
the bicentralizer property of the inclusion of full groups $\Cal N_{M_0}(A_0)\subset \Cal N_{\text{\bf M}}(\text{\bf A})$ 
characterizes the amenability  of $A_0\subset M_0$.

\vskip .05in

$3^\circ$ G. Elek and G. Szabo proved in [ES]  the following ``unique embedding'' type characterization of the amenability property 
for a countable group $H$, analogue to the one for finite separable von Neumann algebras in [J]: if $H$ is amenable then any two embeddings of $H$ into the normalizer $\Cal N$ of 
$\text{\bf A}=\Pi_\omega D_n \subset \Pi_\omega M_{n\times n}(\Bbb C)=\text{\bf M}$, acting freely on $\text{\bf A}$, are conjugate 
by a unitary in $\Cal N$ (this easily implies the same thing for $\text{\bf A}=D^\omega \subset R^\omega =\text{\bf M}$; note that 
by Corollary 5.2 below, the same ``unique embedding'' result actually holds true for ANY ultraproduct inclusion $\text{\bf A}\subset \text{\bf M}$); and that 
if $H$ is sofic and non-amenable, then there exist at least two embeddings of $H$ into $\Cal N$ acting freely on $\text{\bf A}$, non-conjugate 
by unitaries in $\Cal N$. In fact,  as we mentioned in 2.3.2$^\circ$ above, by (8.1 in  [B]) there even exist uncountably many non-conjugate such embeddings. 

Part $1^\circ$ of Theorem 2.2  suggests the following alternative ``bicentralizer'' characterization of amenability for countable groups: 

\vskip .05in 
\noindent 
{\bf (2.3.3)}  {\it Conjecture}:  Let $H$ be a countable group embeddable into the normalizer 
of an ultraproduct MASA $\text{\bf A}\subset \text{\bf M}$ 
(notably $D^\omega \subset R^\omega$, or $\Pi_\omega D_n \subset \Pi_\omega M_{n\times n}(\Bbb C)$),  
such that $H$ acts freely on $\text{\bf A}$ and such that  
it satisfies the bicentralizer condition $(H'\cap \text{\bf A})'\cap \text{\bf M}=\text{\bf A}\vee H$.  
Then $H$ is amenable. 

\heading 3. Approximate free independence in subalgebras 
\endheading

\noindent
\vskip .05in \noindent {\bf 3.1. Notation}. Let $\Cal M$ be a von Neumann algebra.  If $v \in \Cal M$ is a
partial isometry with $v^*v = vv^*$, $X \subset \Cal M$ is a subset and
$k$ a nonnegative integer, then denote $X_v^{0} \overset
\text{\rm def} \to = X$ and $X_v^{k} \overset \text{\rm def} \to =
\{ x_0 \overset k \to{\underset i = 1 \to \Pi} v_i x_i \mid x_i
\in X, \ 1 \leq i \leq k - 1, \ x_0, \ x_k \in X \cup \{ 1 \}, 
v_i\in \{v, v^*\}  \}$.

\proclaim{3.2.  Lemma}  Let $Q\subset M $ be an inclusion of  $\text{\rm II}_1$ von Neumann 
algebras and  assume $Q \not\prec_M Q'\cap M$. Let $f\in Q$ be a non-zero projection. 
For any finite set $F\subset M \ominus (Q'\cap M)$, any $n\geq 1$ and any $\varepsilon > 0$, there exists a partial isometry $v$ in $fQf$  
such that  $vv^*=v^*v$,  $\tau(vv^*) > \tau(f)/4$ and 
$\| E_{Q'\cap M} (x) \|_1 \leq \varepsilon$, $\forall x \in \overset
n \to{\underset k = 1 \to \cup} F_v^{k}$. 
\endproclaim
\noindent 
{\it Proof}.  It is clearly sufficient to prove the statement in case $F=F^*$ and 
$\|x\|\leq 1$, $\forall x\in  F$. Let $\delta > 0$. Denote $\varepsilon_0 = \delta, \
\varepsilon_k = 2^{k+1} \varepsilon_{k-1}, \ k \geq 1$. Denote $
{\mycal W}$ $= \{ v \in fQf \mid vv^*=v^*v \in \Cal P(Q)$, $ \|E_{Q'\cap M} (x)\|_{_1}$ $ \leq \varepsilon_k \tau(v^*v)$, $\forall 1 \leq k \leq n, 
\forall x \in F_v^{k} \}.$  

Endow ${\mycal W}$ with the order $\leq$ in
which $w_1 \leq w_2$ iff $w_1 = w_2 w_1^*w_1$. $({\mycal W} , \leq)$
is then clearly inductively ordered.  Let $v$ be a maximal element
in ${\mycal W}$. Assume $\tau(v^*v) \leq \tau(f)/4$ and denote $p = f - v^*v$.  
Note that this implies $\tau(vv^*)/\tau(p) \leq 1/3$. 

If $w$ is a
partial isometry in $pQp$ with $q=ww^*=w^*w$ and we let $u = v + w$, 
then for $x = x_0 \overset k \to{\underset i = 1 \to \Pi} u_i x_i
\in F_u^{k}$ we have
$$
x = x_0 \Pi_{i=1}^k v_i x_i + \Sigma_\ell \Sigma_i z_{0,i} \Pi_{j=1}^\ell
w_{i_j} z_{j,i}, \tag 1
$$
where the sum is taken over all $\ell = 1, 2, \dots , k$ and all $i
= (i_1, \dots , i_\ell)$, with $1 \leq
i_1 < \cdots < i_\ell \leq k$, and where $w_{i_j} = w$ (resp. $w_{i_j} = w^*$) whenever $v_{i_j} =
v$ (resp. $v_{i_j} = v^*$),  $z_{0,i} = x_0 v_1 x_1 \cdots x_{i_1-1} p$, $ z_{j,i} = p x_{i_j}
v_{i_j+1} \cdots v_{i_{j+1}-1} x_{i_{j+1}-1}p$, for $1 \leq j < \ell$, 
and $z_{\ell,i} = p x_{i_\ell} v_{i_\ell + 1} \cdots v_k x_k$ .

\vskip .05in

By applying $E_{Q'\cap M}$ to the above equation, then taking $\| \ \|_1$ and applying triangle inequality, we then get:

$$
\|E_{Q'\cap M}(x)\|_1  \leq  \|E_{Q'\cap M}(x_0 \Pi_{i=1}^k v_i x_i)\|_1 + \Sigma_\ell \Sigma_i \| E_{Q'\cap M}(z_{0,i} \Pi_{j=1}^\ell
w_{i_j} z_{j,i})\|_1 \tag 1'
$$

Since $v\in \mycal W$, the first term on the right side in $(1')$ is majorized by $\varepsilon_k \tau(vv^*)$, so we are left with estimating the terms 
$z=z_{0,i} \Pi_{j=1}^\ell
w_{i_j} z_{j,i}$ in the double summation on the right hand side, which all have $\ell \geq 1$ number of appearances of $w$ or $w^*$. 

\vskip .05in 
{\it The case $\ell \geq 2$}. Since for $y_1, y_2, y \in M$ 
with $\|y_1\| \leq 1, \|y_2\| \leq 1$ we have $\|E_{Q'\cap M}(y_1 y y_2)\|_1 \leq \|y_1 y y_2\|_1\leq \|y\|_1$, it follows that 
for any $\ell \geq 2$ 
we have: 

$$
\|E_{Q'\cap M}(z)\|_1=\|E_{Q'\cap M} (z_{0,i} w_{i_1} z_{1, i} w_{i_2} z_{2,i} \dots w_{i_\ell}
z_{\ell, i}) \|_{_1} \tag 2 
$$
$$
\leq \| w_{i_1} z_{1,i} w_{i_2} \|_{_1} 
=  \| q z_{1,i} q \|_1 = \| q z_{1,i} q \|_{1,pMp} \tau(p),   
$$
where $\tau_{pMp}=\tau(p)^{-1}\tau_M$ and $\| \ \|_{1,pMp}$ denotes the $L^1$-norm on $pMp$ associated with this trace. 

By applying  Theorem 1.4 to the inclusion $pQp\subset pMp$ (with its trace $\tau_{pMp}$) and to the 
finite set $X\subset pMp$ of all elements of the form $z_{1,i}-E_{(Q'\cap M)p}(z_{1,i}) \in pMp\ominus (Q'\cap M)p$, 
for some $i = (i_1,
\dots , i_\ell)$, $\ell \geq 2$, we obtain that for any 
$ \alpha > 0$,  there exists $ q\in \Cal P(pQp)$  such that
$$
\| q z_{1,i} q - E_{(Q'\cap M)p} (z_{1,i}) q \|_{_{1,pMp}} < \alpha \tau_{pMp} (q). \tag 3
$$ 

Thus, by combining $(2)$ and $(3)$  we get  
$$
\|E_{Q'\cap M}(z)\|_1 \leq  \| q z_{1,i} q \|_{1,pMp} \tau(p) \tag 4 
$$
$$
\leq ( \|E_{(Q'\cap M)p} (z_{1,i}) q \|_{1,pMp} + \alpha \tau_{pMp}(q) ) \tau(p)  
$$
$$
= \| E_{(Q'\cap M)p}(z_{1, i}) \|_{1,pMp} \tau_{pMp}(q)  \tau(p) + \alpha \tau(q) 
$$
$$
=\| E_{(Q'\cap M)p}(z_{1, i}) \|_{1,pMp} \tau(q) +\alpha \tau(q).  
$$

We now take into account that by the definition of the norm $\| \ \|_1$, we have 
$$
\| E_{(Q'\cap M)p}(z_{1, i}) \|_{1,pMp} = \sup \{|\tau(yz_{1,i})|/\tau(p) \mid y\in (Q'\cap M)p, \|y\|\leq 1\} \tag 5 
$$
$$
= \sup \{|\tau(y (1-vv^*) x_{i_1} v_{i_1+1} \cdots v_{i_2-1} x_{i_2-1}(1-vv^*))|/\tau(p) \mid y\in Q'\cap M, \|y\|\leq 1\}. 
$$

But since $y\in Q'\cap M$ commutes with $v, 1-vv^*\in Q$ and $\tau$  is a trace, we actually have 
$\tau(y(1-vv^*)x_{i_1} .... x_{i_2-1}(1-vv^*))=\tau(yx_{i_1} .... x_{i_2-1})-\tau(y v^*x_{i_1}... x_{i_2-1}v)$, so the last term in $(5)$ is further majorized by 
$$
\sup \{|\tau(y x_{i_1} v_{i_1+1} \cdots v_{i_2-1} x_{i_2-1})|/\tau(p) \mid y\in Q'\cap M, \|y\|\leq 1\} \tag 6 
$$
$$
+ \sup \{ |\tau(y v^*x_{i_1} v_{i_1+1} \cdots v_{i_2-1}x_{i_2-1}v)/\tau(p)| \mid  y\in Q'\cap M, \|y\|\leq 1\} 
$$
$$
=(\|E_{Q'\cap M}(x_{i_1} v_{i_1+1} \cdots v_{i_2-1} x_{i_2-1})\|_1
$$
$$
 + \|E_{Q'\cap M}(v^*x_{i_1} v_{i_1+1} \cdots v_{i_2-1}x_{i_2-1}v)\|_1)/\tau(p).
$$

Note at this point that $x_{i_1} v_{i_1+1} \cdots v_{i_2-1} x_{i_2-1}$ lies in  $F^{i_2-i_1-1}_v$ and $v^*x_{i_1} v_{i_1+1}$ $ \cdots v_{i_2-1}x_{i_2-1}v$ 
lies in $F^{i_2-i_1+1}_v$.   
Also, $i_2-i_1 + 1 \leq k$, with the only case when $i_2-i_1 + 1 =k$ corresponding to the case $i_1=1$, $i_2=k$,  $l=2$, i.e., 
to the (single) term $z=x_0w_1(px_1v_2x_2 \cdots  v_{k-1}x_{k-1}p)w_kx_k$ of the double summation in $(1')$.   Thus, 
by combining $(4)$ and $(6)$ and using that $\tau(vv^*)/\tau(p) \leq 1/3$ and choosing $\alpha \leq \delta/3$ (which 
is less than  $(\varepsilon_j - \varepsilon_{j-2})/3$, $\forall j$), for this particular $z$ we get 

$$
\|E_{Q'\cap M}(z)\|_1 \leq \varepsilon_{k-2} (\tau(vv^*)/\tau(p))\tau(q) + \varepsilon_k (\tau(vv^*)/\tau(p))\tau(q) + \alpha \tau(q) \tag 7
$$
$$
\leq (\varepsilon_{k-2}/3 + \varepsilon_k/3 + \alpha)\tau(q) \leq  (2\varepsilon_k/3)\tau(q), 
$$
while for any $z$ with $i_2-i_1+1 \leq k-1$, we get 
$$
\|E_{Q'\cap M}(z)\|_1 \leq \varepsilon_{k-3} (\tau(vv^*)/\tau(p))\tau(q) + \varepsilon_{k-1} (\tau(vv^*)/\tau(p))\tau(q) + \alpha \tau(q) \tag 8
$$
$$
\leq (\varepsilon_{k-3}/3 + \varepsilon_{k-1}/3 + \alpha)\tau(q) \leq (2\varepsilon_{k-1}/3)\tau(q).  
$$

Since $2^{k+1}\varepsilon_{k-1}=\varepsilon_k$ and since there are $\overset k \to{\underset i = 2 \to \Sigma} \left( \matrix k \\ i
\endmatrix \right) = 2^k - k -1$ elements in the double sum in (1) for which
$\ell \geq 2$, of which exactly one has  $i_2 - i_1 +1= k$ and the rest satisfy $i_2-i_1+1 \leq k-1$, by summing up $(7)$ and $(8)$, we get
$$
\Sigma_{\ell \geq 2} \Sigma_i \|z_{0,i} \Pi_{j=1}^\ell w_{i_j} z_{j,i}\|_1 \tag 9
$$
$$
\leq (2^k-k-2) (2\varepsilon_{k-1}/3) \tau(q) + (2\varepsilon_k/3)\tau(q)
$$
$$
=\varepsilon_k\tau(q) - (2k+4)(\varepsilon_{k-1}/3)\tau(q). 
$$

\vskip .05in

{\it The case $\ell =1$}. From the double sum on the right hand side of $(1')$ we will now estimate
the terms with $\ell = 1$.  These are terms which are obtained from $x_0
v_1 x_1 v_2 x_2 \dots v_k x_k$ by replacing exactly one $v_i$ by $w_i$, so
they are of the form $z=z_{0,i} w_{i} z_{1,i} $, where $i=1, 2, ..., k$, 
$z_{0,i}=x_0v_1x_1... v_{i-1}x_{i-1}p$, $z_{1,i}=px_i v_{i+1} ... v_k x_k$ and $w_{i}=w^s$ if $v_i=v^s$, $s\in \{\pm 1\}$ (with the convention that 
$v^{-1}=v^*$, $w^{-1}=w^*$). Note that there are $k$ such terms.  

One should notice at this point that in the above estimates we
only used the fact that $w^*w = ww^* = q\in \Cal P(Q)$ and that it satisfies 
$(3)$ for appropriate $\alpha$.  But we did not use so far the actual form of $w$. We will make the appropriate choice for $w$ now, by making use of   
the condition $Q\not\prec Q'\cap M$.  Indeed, by Theorem 1.5 (2.1 in [P10]), this latter condition implies that for all $ \beta > 0$ 
and all finite sets $Y_1=Y^*_1 \subset M \ominus Q'\cap M$, $Y_2=Y_2^* \subset M$, there 
exists a unitary element $w\in qQq$ such that 

$$
\|E_{Q'\cap M}(y_1w y_2)\|_1 < \beta, \|E_{Q'\cap M}(y_2w y_1)\|_1 < \beta, \forall y_1\in Y_1, y_2 \in Y_2. \tag 10
$$
Note that since $Y_1, Y_2$ are selfadjoint sets, by taking  adjoints 
in $(10)$, from these estimates we also get:  
$$
\|E_{Q'\cap M}(y_2w^* y_1)\|_1 < \beta, \|E_{Q'\cap M}(y_1w^* y_2)\|_1 < \beta, \forall y_1\in Y_1, y_2 \in Y_2. \tag 10'
$$

Denote by $Z$ the set of elements of the form 
$x_0v_1x_1... v_{i-1}x_{i-1}$, or $x_i v_{i+1} ... v_k x_k$, for all possible choices arising from elements in $\overset
n \to{\underset k = 1 \to \cup} F_v^{k}$. 
By applying $(10), (10')$ to $\beta  =\varepsilon_{k-1}\tau(q)/2k$, $n\geq 1$ and $Y_2=Z \cup Z^* \cup \{E_{Q'\cap M}(z) \mid z\in Z\cup Z^* \}$, 
$Y_1=\{y_2-E_{Q'\cap M}(y_2) \mid 
y_2 \in Y_2\}$,  it follows that there exists $w\in \Cal U(qQq)$ such that 

$$ 
\|E_{Q'\cap M}((x_0 v_1 x_1 \dots v_{j-1} x_{j-1} \tag 11 
$$
$$
 - E_{Q'\cap M}(x_0 v_1 x_1 \dots v_{j-1} x_{j-1}))w_j x_j v_{j+1} \dots v_k x_k) \|_{_1} \leq \varepsilon_{k-1} \tau(q)/2k,  
$$
$$
\|E_{Q'\cap M}(E_{Q'\cap M}(x_0 v_1 x_1 \dots v_{j-1} x_{j-1}) w_j (x_j v_{j+1} \dots v_k x_k   \tag 11' 
$$
$$
-E_{Q'\cap M}(px_j v_{j+1} \dots v_k x_k))) \|_{_1} \leq \varepsilon_{k-1} \tau(q)/2k.  
$$

Thus, for each element with $\ell=1$ in the double summation $\Sigma_\ell \Sigma_i z_{0,i} \Pi_{j=1}^\ell
w_{i_j} z_{j,i}$ in $(1)$, i.e., of the form $x_0 v_1 x_1 \dots v_{j-1} x_{j-1} w_j x_j v_{j+1} \dots v_k x_k$, 
we have the estimate: 

$$
\|E_{Q'\cap M}(x_0 v_1 x_1 \dots v_{j-1} x_{j-1} w_j x_j v_{j+1} \dots v_k x_k)\|_1 \tag 12 
$$
$$
\leq 2 \varepsilon_{k-1} \tau(q)/ 2k +  \|E_{Q'\cap M}(x_0 v_1 x_1 \dots v_{j-1} x_{j-1})w_j E_{Q'\cap M}(x_j v_{j+1} \dots v_k x_k)\|_1
$$
$$
\leq \varepsilon_{k-1} \tau(q)/k + \gamma,   
$$
where $\gamma$ is the minimum between 
$$\|E_{Q'\cap M}(x_0 v_1 x_1 \dots v_{j-1} x_{j-1})q\|_1=\tau(q)\|E_{Q'\cap M}(x_0 v_1 x_1 \dots v_{j-1} x_{j-1})\|_1$$ 
and 
$$\|qE_{Q'\cap M}(x_j v_{j+1} \dots v_k x_k)\|_1=\tau(q)\|E_{Q'\cap M}(x_j v_{j+1} \dots v_k x_k)\|_1$$

Both elements $x_0v_1 x_1 \dots v_{j-1} x_{j-1}$,  $x_j v_{j+1} \dots v_k x_k$ belong to some $F_v^{j}$ with $j\leq k-1$, and at least one of them 
with $j\neq 0$. Thus, by the properties of $v\in \mycal W$ and the assumption $\tau(vv^*)\leq \tau(f)/4$, 
we have $\gamma \leq \varepsilon_{k-1} \tau(vv^*) \tau(q)\leq \varepsilon_{k-1}\tau(q)/4$. 
Hence, the last term in $(12)$ is majorized by 

$$
\varepsilon_{k-1}\tau(q)/k + 
\varepsilon_{k-1}\tau(q)/4 = (k/4 +1) \varepsilon_{k-1}\tau(q). \tag 13 
$$

\vskip .05in 

{\it Summing up the cases $\ell\geq 2$ and $\ell=1$}. Since there are $k$ terms with $\ell=1$, obtained by taking $j=1, ..., k$, by summing up over $j$ in $(12)-(13)$ 
and combining with the estimate $(9)$, obtained in the case $\ell \geq 2$, we deduce from  $(1')$ the following final estimate: 

$$
\|E_{Q'\cap M}(x)\|_1 \leq  \|E_{Q'\cap M}(x_0 \Pi_{i=1}^k v_i x_i)\|_1  + \Sigma_\ell \Sigma_i \|E_{Q'\cap M}(z_{0,i} \Pi_{j=1}^\ell
w_{i_j} z_{j,i})\|_1 \tag 14
$$
$$
\leq \varepsilon_k \tau(vv^*)+ (\varepsilon_k - (2k+4)\varepsilon_{k-1}/3)  \tau(q) + (k/4 + 1) \varepsilon_{k-1} \tau(q)
$$
$$
\leq \varepsilon_k \tau(vv^*) + \varepsilon_k\tau(ww^*) = \varepsilon_k \tau((v+w)(v+w)^*). 
$$

Since $u=v+w$ has also the property that $uu^*=u^*u$, it follows from $(13)$ that $u\in \mycal W$. 
But this contradicts the maximality of $v \in {\mycal W}$.

We conclude that $\tau(v^*v) >\tau(f)/4$.  If we now take $\delta\leq \varepsilon/2^{n^2+1}$, 
then $\varepsilon_n = 2^{(n+1)(n+2)/2}\delta $ $< 2^{n^2+1} \delta \leq \varepsilon$ and the statement follows.  

\hfill $\square$

\heading 4. Free independence in ultraproduct framework \endheading

\vskip .05in \noindent {\bf
4.1. Notation}. Let $M_n$ be a sequence of finite factors with $\dim(M_n) \rightarrow \infty$. Let  
$\omega$ be a free ultrafilter on $\Bbb N$ and denote $\text{\bf M}=\Pi_\omega M_n$.
We consider the following two special classes of subalgebras  of $\text{\bf M}$: 
\vskip .05in 
\noindent 
$(4.1.1)$ We denote by $\mycal Q_{u}$ the class of von Neumann subalgebras $\text{\bf Q}\subset \text{\bf M}$ 
which are of the form $\text{\bf Q}=\Pi_\omega Q_n$, for some subalgebras $Q_n \subset M_n$, 
and have the property that $\text{\bf Q}\not\prec_{\text{\bf M}} \text{\bf Q}'\cap \text{\bf M}$.  
\vskip .05in 
\noindent
$(4.1.2)$ We denote by $\mycal Q_{b}$ the class of von Neumann subalgebras $\text{\bf Q}\subset \text{\bf M}$ with the property that 
$\text{\bf Q}'\cap \text{\bf M}$ is separable and  $(\text{\bf Q}'\cap \text{\bf M})'\cap \text{\bf M}=\text{\bf Q}$.   

\vskip .1in 

The next result provides some properties and examples of algebras in these two classes. 

\proclaim{4.2. Proposition} $1^\circ$ If $\text{\bf Q}\in \mycal Q_u$, then $\text{\bf Q}$ is of type $\text{\bf II}_1$.

$2^\circ$  If $Q_n \subset M_n$ are von Neumann subalgebras such 
that $Q_n \not\prec_{M_n} Q_n'\cap M_n$, $\forall n$, then $\text{\bf Q}=\Pi_\omega Q_n$ satisfies $\text{\bf Q}\not\prec_{\text{\bf M}} 
\text{\bf Q}'\cap \text{\bf M}$, and thus $\text{\bf Q} \in \mycal Q_{u}$. 

$3^\circ$ Assume $m_n$ is an increasing sequence of positive integers of the form $m_n=d_n \cdot k_n$, with $d_n, k_n \in \Bbb N$.  
Let $M_n=M_{m_n \times m_n}(\Bbb C)$, with $P_n=M_{d_n \times d_n}(\Bbb C)$, $Q_n=M_{k_n \times k_n}(\Bbb C)$, viewed as 
subalgebras of $M_n$ that commute and generate $M_n$. Then $\text{\bf Q}=\Pi_\omega Q_n$, $\text{\bf P}=\Pi_\omega P_n$  
satisfy the following properties: $\text{\bf Q}'\cap \text{\bf M}=\text{\bf P}$, $\text{\bf P}'\cap \text{\bf M}=\text{\bf Q}$;  
$\text{\bf Q} \not\prec_{\text{\bf M}} \text{\bf P}$ $($and thus $\text{\bf Q} \in \mycal Q_{u})$ if and only if $\lim_\omega d_n/k_n =  0$. 

$4^\circ$ If $B\subset \text{\bf M}$ is a separable amenable von Neumann subalgebra, then $\text{\bf Q}:=B'\cap \text{\bf M}$ 
satisfies $\text{\bf Q}'\cap \text{\bf M}=B$. Thus $\text{\bf Q} \in \mycal Q_b$. 

$5^\circ$ If $\text{\bf Q} \in \mycal Q_b$ then $\text{\bf Q}$ is of type $\text{\bf II}_1$, has no separable direct summand, and $\text{\bf Q}\not\prec_{\text{\bf M}} \text{\bf Q}'\cap 
\text{\bf M}$ $($the latter being separable$)$. 
\endproclaim 
\noindent
{\it Proof}. $1^\circ$ If an inclusion of finite von Neumann algebras   $B\subset M$ is so that $B$ is type I, then 
there exists a non-zero projection $e\in B$ such that $eBe$ is abelian, implying that $eBe\subset (eBe)'\cap eMe$, 
thus $B \prec_M B'\cap M$. Since in our case we have $\text{\bf Q}\not\prec_{\text{\bf M}} \text{\bf Q}'\cap \text{\bf M}$, this shows 
that $\text{\bf Q}$ cannot have type I summands, thus $\text{\bf Q}$ is  II$_1$.  

Part $2^\circ$ is an immediate consequence of Theorem 1.5 and of the fact that $\text{\bf Q}'\cap \text{\bf M}=\Pi_\omega (Q_n'\cap M_n)$ 
with $E_{\text{\bf Q}'\cap \text{\bf M}}(x)=(E_{Q_n'\cap M_n}(x_n))_n$, for $x=(x_n)_n\in \text{\bf M}=\Pi_\omega M_n$. 

Part $3^\circ$ is an easy exercise (using Theorem 1.5) while part $4^\circ$ is a direct consequence of Theorem 2.1. 

To prove part  $5^\circ$, note that if $\text{\bf Q} \in \mycal Q_b$ then $\text{\bf Q}$  has no separable direct summand, by the same 
observation we have used in the proof of part $1^\circ$.  

\hfill 
$\square$

\vskip .05in 
Note that conjecture $(2.3.1)$ predicts that  in fact the class $\mycal Q_b$ only consists of centralizers of separable amenable subalgebras of $\text{\bf M}$,   
i.e., that any subalgebra in $\mycal Q_b$ is of the form $4.2.4^\circ$ above.   

\proclaim{4.3. Theorem} Assume $\text{\bf Q} \subset \text{\bf M}$ is either in the class $\mycal Q_u$, or $\mycal Q_b$. 
If $X\subset \text{\bf M}\ominus (\text{\bf Q}'\cap \text{\bf M})$ is a separable subspace, then there exists a diffuse abelian von Neumann subalgebra 
$A\subset \text{\bf Q}$  such that $A$ is free independent to $X$, relative to $\text{\bf Q}'\cap \text{\bf M}$, more precisely  
$E_{\text{\bf Q}'\cap \text{\bf M}}(x_0 \overset n \to{\underset i = 1 \to \Pi} a_ix_i )=0$, for all $n \geq 1$ and all 
$x_i \in X$, $1\leq i \leq n-1$, $x_0, x_n \in X \cup \{1\}$, 
$a_i \in A\ominus \Bbb C$, $1\leq i \leq n$. 
\endproclaim

\proclaim{4.4. Corollary} With the same assumptions and notations as in $4.3$ above, we have: 

\vskip .05in
$1^\circ$ Let $P\subset \text{\bf M}$ be a von Neumann subalgebra making a commuting square with $\text{\bf Q}'\cap \text{\bf M}$ 
and denote $B_1=P \cap (\text{\bf Q}'\cap \text{\bf M})$. Assume that $L^2P$ is countably generated both as a left and as a right $B_1$ Hilbert module 
$($equivalently, there exists a separable space $X \subset P$ such that $X\perp B_1$,  
and $\text{\rm sp} X B_1$ and $\text{\rm sp} B_1 X$ are both $\| \ \|_2$-dense in $P\ominus B_1)$. 
Then there exists a diffuse von Neumann subalgebra $B_0\subset \text{\bf Q}$ such that $P\vee B_0 \simeq P *_{B_1} (B_1 \overline{\otimes} B_0)$. 
\vskip .05in 
$2^\circ$ Let $N_i\subset \text{\bf M}$ be separable von Neumann algebras, with amenable subalgebras $B_i$, $i=1,2$, such that $(B_1, \tau)\simeq (B_2,\tau)$. 
Then there exists  a unitary element $u\in \text{\bf M}$ such that $uB_1u^*=B_2$ and such that, after identifying $B=B_1\simeq B_2$ via $\text{\rm Ad}(u)$, 
we have $N_1 \vee uN_2u^* \simeq N_1*_B N_2$. 
\endproclaim

Note that the case $B$ atomic of 4.4.2$^\circ$ above has already been shown in [P6]. The case  
$\text{\bf M}=R^\omega$ of 4.4.2$^\circ$ shows in particular that if $N_1, N_2$ are two separable finite von Neumann algebras 
with a common amenable subalgebra $B\subset N_i$ and $N_1, N_2$ are both embeddable into $R^\omega$, then so is $N_1 *_B N_2$.  
This recovers a result in  [BDJ]  (see also [FGR] for more recent related considerations). 

A particular case when the assumptions in $4.4.1^\circ$ are satisfied, is when the subalgebra 
$P\subset \text{\bf M}$ making a commuting square 
with $\text{\bf Q}'\cap \text{\bf M}$ is itself separable. But there are interesting non-separable examples as well, 
that can even allow obtaining free product with amalgamation over the entire $\text{\bf Q}'\cap \text{\bf M}$ (which is non-separable in case $\text{\bf Q}\in\mycal Q_u$). 
For instance, if $\Cal U\subset \Cal U(\text{\bf M})$ is a countable group of unitaries normalizing $\text{\bf Q}'\cap \text{\bf M}$, then the von Neumann 
algebra $P$ generated by $\Cal U$ and $\text{\bf Q}'\cap \text{\bf M}$ satisfies all the conditions in $4.4.1^\circ$ with $B_1=\text{\bf Q}'\cap \text{\bf M}$. 

Note in this respect that one can alternatively  take in the statement of Theorem 4.3 the separable space $X$ to be of the form 
$X=P\ominus (P\cap \text{\bf Q}'\cap \text{\bf M})$, for some 
separable von Neumann algebra $P$ making a commuting square with $\text{\bf Q}'\cap \text{\bf M}$. Indeed, due to Lemma 1.2, the two versions 
follow equivalent.

\proclaim{Lemma 4.5} Let $\text{\bf Q}\subset \text{\bf M}$ be a von Neumann subalgebra lying in 
either the class $\mycal Q_u$ or the class $\mycal Q_b$. Let $f\in \text{\bf Q}$ be a non-zero projection 
and $X\subset \text{\bf M}\ominus \text{\bf Q}'\cap \text{\bf M}$ 
a countable set. 
Then there exists a partial isometry $v$ in $f\text{\bf Q}f$ such that $vv^*=v^*v$, 
$\tau(vv^*) \geq \tau(f)/4$ and $E_{\text{\bf Q}'\cap \text{\bf M}}(x)=0$,  $\forall x \in X_v^{k}$, $\forall k \geq 1$. 
\endproclaim 
\noindent
{\it Proof}. Let $X = \{ x_k \}_{k\geq 1}$ be an enumeration of $X$ and denote $x_0=1$. 
By applying Lemma 3.2 to the inclusion of II$_1$ von Neumann algebras $\text{\bf Q} \subset
\text{\bf M}$, the projection $f\in \text{\bf Q}$, the positive constant $\varepsilon = 2^{-n}$ and the finite set $X_n = \{ x_k \mid k \leq
n\}$,  we get a partial isometry $w_n$ in $f\text{\bf Q}f$ with the property that $w_nw_n^*=w_n^*w_n$, 
$\tau(w_n^* w_n)\geq \tau(f)/4$  and 

$$
\| E_{\text{\bf Q}'\cap \text{\bf M}} (x) \|_1 < 2^{-n}, \forall x \in \underset
k \leq n \to \cup (X_n)_{w_n}^{k}. \tag 1 
$$

Let  $f=(f_m)_m$ 
be a representation of $f$ with $f_m$ projections. Let also $x_k=(x_{k,m})_m$ be a representation 
of $x_k$, with $x_{k,m}\in M_m$, $\|x_{k,m}\|\leq \|x_k\|$, $\forall k, m$, and $w_k=(w_{k,m})_m\in \text{\bf Q}$ a representation of $w_k$ 
with $w_{k,m}$  partial isometries satisfying $w_{k,m}w_{k,m}^*=w_{k,m}^*w_{k,m} \leq f_m$.

Assume first that $\text{\bf Q}=\Pi_\omega Q_n \in \mycal Q_u$, 
in which case we may clearly also assume $f_m \in  \Cal P(Q_m)$ and $w_{k,m}\in f_mQ_mf_m$, $\forall k, m$. Noticing 
that if $y=(y_n)_n \in \text{\bf M}$ then $E_{\text{\bf Q}'\cap \text{\bf M}}(y)=(E_{Q_n'\cap M_n}(y_n))_n$,   
it follows from $(1)$ that 

$$
\underset m \rightarrow \omega \to \lim \|E_{Q_m'\cap M_m}(x_{j_0,m}\Pi_{i=1}^k w_{n,m}^{\gamma_i} x_{j_i,m})\|_1 < 2^{-n}, \tag 2
$$
for all $1\leq k \leq n$, $x_{j_0}, x_{j_k}\in X_n \cup \{1\}$, $x_{j_i}\in X_n$, $1\leq i \leq k-1$, $\gamma_i\in \{\pm 1\}$, $\forall i$ (as before,  
for partial isometries $y\in \text{\rm M}$ with $yy^*=y^*y$, we use the convention $y^{-1}=y^*$).

Let $V_n$ be the set of all $m\in \Bbb N$ with the property that 

$$
\|E_{Q_m'\cap M_m}(x_{j_0,m}\Pi_{i=1}^k w_{n,m}^{\gamma_i} x_{j_i,m})\|_1 < 2^{-n}, \tag 3 
$$
for all $1\leq k \leq n$, $1 \leq j_i \leq n$ for $i\geq 1$, $0\leq j_0\leq n$, $\gamma_i \in \{\pm 1\}$. By $(2)$ it follows that 
$V_n$ corresponds to an open-closed neighborhood of $\omega$ in $\Omega$, 
under the identification $\ell^\infty \Bbb N = C(\Omega)$. Let now $W_n$, $n\geq 0$,  be defined 
recursively as follows: $W_0=\Bbb N$ and $W_{n+1}=W_n \cap V_{n+1} \cap \{n \in \Bbb N \mid n > \min  W_n\}$. Note that, 
with the same identification as before, $W_n$ is 
a strictly decreasing sequence of neighborhoods of $\omega$. 

Define $v=(v_m)_m$ by letting $v_m=w_{n,m}$ for $m\in W_{n-1} \setminus W_n$. 
It is then easy to check that  $v$ is a partial isometry in $f\text{\bf Q}f$ satisfying all  the required conditions.  

Assume now that $\text{\bf Q}\in \mycal Q_b$. Let $\{y_\ell\}_\ell \subset \text{\bf Q}'\cap \text{\bf M}$ be a countable set dense 
in the unit ball of $\text{\rm Q}'\cap \text{\bf M}$ in the norm $\| \ \|_2$. 
Note that if $y_\ell=(y_{\ell,m})_m$ then $x=(x_n)_n \in \text{\bf M}$ satisfies $x\in \text{\bf Q}$ 
iff $\underset m \rightarrow \omega \to \lim \|[x_m, y_{\ell,m}]\|_2=0$, $\forall \ell$. 
Also, $x \perp \text{\bf Q}'\cap \text{\bf M}$ iff $\underset m \rightarrow \omega \to \lim \tau(x_my_{\ell,m})=0$, $\forall \ell$. Moreover, if $\delta > 0$, then 
$\|E_{\text{\bf Q}'\cap \text{\bf M}}(x)\|_1\leq \delta$ iff $\underset m \rightarrow \omega \to \lim |\tau(x_my_{\ell,m})| \leq \delta$, $\forall \ell$. 

With this in mind, from $(1)$ it follows that  the partial isometries $w_n=(w_{n,m})_m \in \text{\bf Q}$ satisfy

$$
\underset m \rightarrow \omega \to \lim |\tau((x_{j_0,m}\Pi_{i=1}^k w_{n,m}^{\gamma_i} x_{j_i,m})y_{\ell,m})| < 2^{-n}, \tag 4
$$
for all $1\leq k \leq n$, $x_{j_0}, x_{j_k} \in X_n \cup \{1\}$, $x_{j_i}\in X_n$, $1\leq i \leq k-1$, $\gamma_i\in \{\pm 1\}, \forall i$, and for all $\ell \geq 1$.
Also, the fact that $w_n$ belongs to $f\text{\bf Q}f$ is equivalent to 
$$
\underset m \rightarrow \omega \to \lim \|[w_{n,m}, y_{\ell,m}]\|_2=0, \forall \ell; \underset m \rightarrow \omega \to \lim \|f_mw_{n,m}f_m - w_{n,m}\|_1=0 \tag 5
$$

Let $V_n$ be the neighborhood of $\omega$ consisting of all $m\in \Bbb N$ with the property that 

$$
|\tau((x_{j_0,m}\Pi_{i=1}^k w_{n,m}^{\gamma_i} x_{j_i,m})y_{\ell,m})| < 2^{-n}; \tag 6
$$ 
$$
 \|[w_{n,m}, y_{\ell,m}]\|_2 < 2^{-n}; \|f_mw_{n,m}f_m - w_{n,m}\|_1< 2^{-n};  
 $$
 for all $\ell = 1, 2, ..., n$ as well as for all $1\leq k \leq n$, $x_{j_0}\in X_n \cup \{1\}$, $x_{j_i}\in X_n$, $\gamma_i\in \{\pm 1\}$. Let further 
 $W_n\subset \Bbb N$, $n\geq 0$, be defined recursively as follows: $W_0=\Bbb N$ and $W_{n+1} = W_n\cap V_{n+1} 
 \cap \{n \in \Bbb N \mid n > \min W_n\}$. It follows  that $W_n$ are all neighborhoods of $\omega$, that $W_n \subset \cap_{j\leq n} V_j$, $W_{n+1}\subset W_n$, 
 and $W_{n+1}\neq W_n$. 
  
We now define $v=(v_m)_m$, by letting $v_m=w_{n,m}$ if $m\in W_{n-1} \setminus W_n$. By the way $w_{n,m}$ have been taken,     
 $v$ follows a partial isometry with $vv^*=v^*v$, while by the second relation in $(6)$ we have 
 $v\in f\text{\bf Q}f$ and by the first relation  in $(6)$ 
 we have $E_{\text{\bf Q}'\cap \text{\bf M}}(x)=0$,  $\forall x \in X_v^{k}$, $\forall k \geq 1$. 
\hfill
$\square$

\vskip .05in
\noindent {\it Proof of} 4.3. Recall that by 4.2.1$^\circ$ and 4.2.5$^\circ$, $\text{\bf Q}$ is of type II$_1$. Thus it contains a copy $R$ of 
the hyperfinite II$_1$ factor $R$ and any element $y \in R$ of trace $0$ satisfies $E_{\text{\bf Q}'\cap \text{\bf M}}(y)=0$. Fix 
$u_0\in R$ to be a Haar unitary and denote $Y=X \cup \{u_0^k \mid k\neq 0\}$. Note that $Y$ is a countable subset   
in $\text{\bf M}\ominus (\text{\bf Q}'\cap \text{\bf M})$. We construct recursively a sequence 
of partial isometries $v_1, v_2, .... \in \text{\bf Q}$ such that 
\vskip .05in 
\noindent
$(i)$ $v_{j+1}v_j^*v_j = v_j$, $v_jv_j^*=v_j^*v_j$ and $\tau (v_jv_j^*) \geq 1 - 1/2^j$, $\forall j \geq 1$. 

 \vskip .05in 
\noindent
$(ii)$ $E_{\text{\bf Q}'\cap \text{\bf M}}(x)=0$,  
$\forall x \in Y_{v_j}^{k}$, $\forall k\geq 1$. 

\vskip .05in 

Assume we have constructed $v_j$ for $j=1, ..., m$. If $v_m$ is a unitary element, then we let $v_j=v_m$ for all $j\geq m$. 
If $v_m$ is not a
unitary element, then let $f = 1 - v_m^*v_m \in \text{\bf Q}$. Note that 
$E_{\text{\bf Q}'\cap \text{\bf M}} (x') = 0$, for all $x' \in X' \overset \text{\rm def}
\to =\cup_k Y_{v_m}^{k}$. Thus, if we apply Lemma 4.5 to $\text{\bf Q} \subset \text{\bf M}$, the projection $f\in \text{\bf Q}$ and the countable set 
$X' \subset \text{\bf M}\ominus (\text{\bf Q}'\cap \text{\bf M})$, then we get a partial isometry $u \in f\text{\bf Q}f$, with $uu^*=u^*u$ 
satisfying $\tau(uu^*)\geq \tau(f)/4$ and $E_{\text{\bf Q}'\cap \text{\bf M}} (x) = 0$ for all $x \in 
\cup_k (X')_u^{k}$.  But then $v_{m+1}=v_m + u$ will satisfy  both $(i)$ and $(ii)$ for $j=m$. 

It follows now from $(i)$ that the sequence $v_j$ converges in the norm $\|$ $\|_2$ to a unitary element $v\in \text{\bf Q}$, which due to $(ii)$ will satisfy the condition  
$E_{\text{\bf Q}'\cap \text{\bf M}}(x)=0$, $\forall x \in \cup_n Y^n_v$. But then the von Neumann algebra $A$ 
generated by the Haar unitary $u=vu_0v^*\in \text{\bf Q}$ clearly satisfies 
the conditions required in $4.3$.  
\hfill $\square$

\vskip .05in
\noindent
{\it Proof of} 4.4. $1^\circ$ Let $X_0\subset P\ominus B_1$ be a separable subspace such that sp$X_0B_1$ and sp$B_1X_0$  are $\| \ \|_2$-dense 
in $P\ominus B_1$. By Theorem 4.3, there exists a diffuse von Neumann subalgebra $B_0\subset \text{\bf Q}$ such that $B_0$ is free 
independent to $X_0$ relative to $\text{\bf Q}'\cap \text{\bf M}$. It is sufficient  to prove that $E_{\text{\bf Q}'\cap \text{\bf M}}(x_0\Pi_i y_ix_i) =0$, 
for any $x_0\in X_0B_1 \cup \{1\}$, $x_i \in X_0B_1, y_i \in B_0\ominus \Bbb C1$, $1\leq  i \leq n$. But any element in $X_0B_1$ can be approximated arbitrarily well 
by a linear combination of elements in $B_1X_0$. The ``coefficient'' in $B_1$ of each one of these elements 
commutes with $y_{i-1}$, so we can ``move it to the left'', being ``swollen'' by the $x_{i_1}\in X_0B_1$. Thus, in the end, it follows that it is sufficient to have 
$E_{\text{\bf Q}'\cap \text{\bf M}}(x_{0,0} \Pi_i y_i x_{0,i})=0$ for $x_{0,0}\in X_0 \cup \{1\}$, $x_{0,i}\in X_0$, $y_i \in B_0\ominus \Bbb C1$, which is indeed the case 
because $B_0$ is free independent to $X_0$ relative to $\text{\bf Q}'\cap \text{\bf M}$.

$2^\circ$ By the first part of Remark 2.3, after possibly conjugating with a unitary $u_0\in \text{\bf M}$, we 
may assume the subalgebras $B_1$, $B_2$ coincide. Denote $B$ this common algebra and let $\text{\bf Q}=B'\cap \text{\bf M}$, which by 2.1 satisfies $\text{\bf Q}'\cap 
\text{\bf M}=B$ and by 4.2.4$^\circ$ it belongs to $\mycal Q_b$. Now apply 
$4.3$ to $\text{\bf Q}$ and to the separable space $X=N_1 \ominus B + N_2\ominus B$, to conclude that there exists a unitary element $u\in \text{\bf Q}$ such that  
$uN_2u^*$ and $N_1$ generate the free amalgamated product $\simeq N_1*_B N_2$. 
\hfill
$\square$

\heading 5.   More on the incremental patching method \endheading

The crucial step in proving Theorem 4.3 is Lemma 3.2. The technique used in its proof 
consists of building unitaries $u$ that are approximately $n$-independent with respect to certain finite sets,  
by ``patching'' together infinitesimal pieces of $u$. This technique was  
first considered  in (2.1 of [P3]), to show that given any countable set $X$ in a finite von Neumann algebra $M$ 
and any diffuse abelian von Neumann subalgebra $A \subset M$, there exists a Haar unitary $u\in A^\omega$ such that any word 
that alternates letters from $X$ and $\{u^n \mid n \geq 1 \}$, has $0$-trace. This result was a key tool in proving 
that any derivation of a II$_1$ factor into the ideal of compact operators is inner, in [P3].   

The technique was substantially refined in [P6], to prove a particular case of the case $Q\in \mycal Q_u$ of   
Theorem 4.3, in which $\text{\bf Q}=\Pi_\omega Q_n \in \mycal Q_u$ 
is so that $Q_n \subset M_n$ are II$_1$ subfactors with atomic relative commutant $Q_n'\cap M_n$ (which thus clearly satisfy $Q_n\not\prec_{M_n} Q_n'\cap M_n$). 
The result in [P6] had several applications over the years. For instance, it played an important role in 
developing reconstruction methods in Jones theory of subfactors in ([P4], [P7], [P9]) and it 
led, in combination with ([V]),  to the definition of amalgamated free product  
of inclusions of finite von Neumann algebras in [P4]. 
It was also used to prove key technical results in ([IPP], [Va]) and to 
show that the free product of standard invariants of subfactors defined in ([BiJ]) can be realized in the hyperfinite II$_1$ factor $R$ 
(see A.2 in [IPP] and [Va]).  

More recently, the same incremental patching method was used in [P12] to prove that if $A_n \subset M_n$ 
is a sequence of MASAs in II$_1$ factors, then 
the abelian von Neumann algebra $\text{\bf A}=\Pi_\omega A_n \subset \Pi_\omega M_n=\text{\bf M}$ 
contains diffuse subalgebras $B_0$ that are $\tau$-independent to any given separable subalgebra $B\subset \text{\bf A}$ 
and $3$-independent to any given countable set $X\subset \text{\bf M}\ominus \text{\bf A}$, 
i.e. any alternating word with at most 3 letters from $X$ and 3 letters from $B_0\ominus \Bbb C1$ has trace $0$ (see 0.2 in [P12]). Moreover, 
if $A_n$ are all {\it singular} (in the sense of [D1], i.e. any unitary normalizing $A_n$ is contained in $A_n$), 
then $B_0$ can be chosen to be free independent to $X$, relative to $\text{\bf A}$, a fact that allowed settling the Kadison-Singer 
problem for ultraproducts of singular MASAs $\text{\bf A}\subset \text{\bf M}$ (see 0.1 in [P12]). 

A concrete example of a diffuse subalgebra $B_0$ in an ultraproduct MASA $\text{\bf A}$  
satisfying the 3-independence property is the following: 

Let $\Gamma \curvearrowright X$ be an ergodic (but not necessarily free)  
measure preserving action of a discrete group $\Gamma$ on a probability space $(X, \mu)$ and $\Gamma \curvearrowright Y=[0, 1]^\Gamma$ 
be the Bernoulli $\Gamma$-action with diffuse base. Let $A=L^\infty(X)\otimes L^\infty(Y)$ with $\Gamma \curvearrowright A$ the product action. 
Let $M = A \rtimes \Gamma$ and $\text{\bf A}=A^\omega \subset M^\omega = \text{\bf M}$. 
If we take $B=L^\infty(X)$ and let $B_0=1 \otimes L^\infty([0,1]) \otimes 1 \subset L^\infty(Y)$ be 
the base of the Bernoulli action, viewed as a tensor component of the infinite tensor product $L^\infty(Y)= \otimes_{g\in \Gamma} (L^\infty([0,1]))_g$, then  
it is easy to see that $B_0$ is $\tau$-independent to $B$ and $3$-independent with respect to $X=\{u_g \mid g\in \Gamma\}$.  

This construction can actually be recovered ``asymptotically'' 
inside any group measure space von Neumann algebra. Indeed, using the incremental patching  technique, we will now  
prove that  (generalized) Bernoulli $\Gamma$-actions can be retrieved 
inside any free action of $\Gamma$ on an ultrapower of measure spaces. More generally  we have:

\proclaim{5.1. Theorem} Let $A_n \subset M_n$ be a sequence of MASAs in finite factors, with $\text{\rm dim} M_n \rightarrow \infty$, and  
denote $\text{\bf A}=\Pi_\omega A_n \subset \Pi_\omega M_n = \text{\bf M}$. Assume $\Gamma \subset \Cal N_{\text{\bf M}}(\text{\bf A})$ 
is a countable group of unitaries acting freely on $\text{\bf A}$ and let $H\subset \Gamma$ be an amenable subgroup. Given any separable abelian von Neumann 
subalgebra $B \subset \text{\bf A}$, there exists a separable diffuse abelian subalgebra $A\subset \text{\bf A}$ such that:  
$A, B$ are $\tau$-independent, $\Gamma$ normalizes $A$, and the action of $\Gamma$ on $A$ is isomorphic to the generalized Bernoulli 
action $\Gamma \curvearrowright L^\infty([0,1]^{\Gamma/H})$. 
\endproclaim 
\noindent
{\it Proof}. Let $\{u_g \mid g\in \Gamma \}$ be the unitaries in $\Gamma$. Denote by $g_0=1, g_1, g_2, ... \in \Gamma$ 
a set of representants of $\Gamma/H$.  It is clearly sufficient to construct a Haar unitary $w$ in $\text{\bf A}$ such that 
$w$ commutes with $u_h, \forall h\in H$, and such that $B$ and 
$u_{g_i} \{w^n \mid n\in \Bbb Z\} u_{g_i}^*$, $i=0, 1, 2, ...$, are all multi-independent, in the sense that for any $k$, 
any non-zero integers $n_j$, distinct non-negative integers $m_j$, and any $b\in B$, we have $\tau(b \Pi_{j=0}^k u_{g_{m_j}}w^{n_j}u_{g_{m_j}}^*)=0$. 

We need some notations. Thus, we let $\text{\bf A}_0$ be the subalgebra of all elements in $\text{\bf A}$ 
that are fixed by $H$ and let  $\{b_n\}_n$ be a $\| \ \|_2$-dense subset of 
the unit ball of $B$. If $v$ is a partial isometry in $\text{\bf A}_0$, then we denote by $F_{v,n}$ the set of all elements of the form 
$b_i \Pi_{j=0}^k u_{g_{m_j}}v^{n_j}u_{g_{m_j}}^*$, 
where $1\leq i\leq n$, $1\leq k \leq n$, $m_j$ are distinct integers beween $0$ and $n$, and $1\leq |n_j|\leq n$. We first prove the following: 

\vskip .05in
{\it Fact}. Given any $n\geq 1$ and any $\delta > 0$, there exists a 
Haar unitary $v\in \text{\bf A}_0$ such that  $|\tau(x)|\leq \delta$, $\forall x\in F_{v,n}$. 

\vskip .05in 
To prove this, let ${\mycal W}:=\{ v \in \text{\bf A}_0 \mid  |\tau(x)| \leq \delta \tau(v^*v), 
\forall x \in F_{v,n}, \tau(v^m)=0, \forall m\neq 0 \}$.  Endow ${\mycal W}$ with the order $\leq$ in
which $w_1 \leq w_2$ iff $w_1 = w_2 w_1^*w_1$. $({\mycal W} , \leq)$
is then clearly inductively ordered.  Let $v$ be a maximal element
in ${\mycal W}$. Assume $\tau(v^*v) < 1$ and denote $p = 1 - v^*v$. If $w\in \text{\bf A}_0p$ is a 
partial isometry satisfying $ww^*=w^*w$, $\tau(w^m)=0$, $\forall m\neq 0$, and we denote $u=v+w$, 
then by noticing that $(v+w)^{n_j}=v^{n_j}+w^{n_j}$, we obtain:  
$$
b_i \Pi_{j=0}^k u_{g_{m_j}}u^{n_j}u_{g_j}^*=b_i \Pi_{j=0}^k u_{g_{m_j}}v^{n_j}u_{g_{m_j}}^*+ \Sigma b_i  \Pi_{j=0}^k u_{g_{m_j}}z_j^{n_j}u_{g_{m_j}}^*,  \tag 1
$$
where $z_j\in \{v, w\}$ and the sum is taken over all possible choices for $z_j=v$ or $z_j=w$, with at least one occurrence of $z_j=w$ (thus, there are  $2^{k+1}-1$ 
many terms in the summation). We thus get the estimate

$$
|\tau(b_i \Pi_{j=0}^k u_{g_{m_j}}u^{n_j}u_{g_{m_j}}^*)| \tag 2
$$
$$
\leq |\tau(b_i \Pi_{j=0}^k u_{g_{m_j}}v^{n_j}u_{g_{m_j}}^*)|+ 
\Sigma|\tau(b_i  \Pi_{j=0}^k u_{g_{m_j}}z_j^{n_j}u_{g_{m_j}}^*)|  
$$
$$
\leq \delta \tau(vv^*) + \Sigma' |\tau(b_i  \Pi_{j=0}^k u_{g_{m_j}}z_j^{n_j}u_{g_{m_j}}^*)| +
 \Sigma'' |\tau(b_i   \Pi_{j=0}^k u_{g_{m_j}}z_j^{n_j}u_{g_{m_j}}^*)| 
$$
where the summation $\Sigma'$ contains the terms with just one occurrence of  $z_j=w$ and $\Sigma''$ is the summation of the terms that have at least 
2 occurrences of $z_j=w$. Since $\text{\bf A}$ is abelian, the terms $u_{g_{m_j}}z_j^{n_j}u_{g_{m_j}}^*$ in a product can be permuted arbitrarily. Thus, in 
each summand of $\Sigma''$ we can bring two of the occurrences of $w$ so that to be adjacent, 
i.e., of the form $y_1u_{g_{m_j}}w^{n_j}u_{g_{m_j}}^*u_{g_{m_i}}w^{n_i}u_{g_{m_i}}^*y_2$. Recall at this point that by Theorem 2.2.1$^\circ$ we have 
$\text{\bf A}_0'\cap \text{\bf M}= \text{\bf A}\vee H$. 
Since $g_{m_i} \neq g_{m_j}$ for all $i\neq j$, by applying part $1^\circ$ of Theorem 1.4 to $Q=\text{\bf A}_0p$ and the finite set $F=\{u_{g_{m_j}}^*u_{g_{m_i}} \mid i\neq j\} 
\perp \text{\bf A}\vee H = \text{\bf A}_0'\cap \text{\bf M}$, it follows that for any $\alpha > 0$, 
there exists a non-zero  $q\in \Cal P(\text{\bf A}_0p)$ such that 

$$
\|qu_{g_{m_j}}^*u_{g_{m_i}}q\|_1 < \alpha \tau(q), \forall 0 \leq m_i \neq m_j \leq n. \tag 3 
$$
 
Since there are $2^{k+1}-(k+1)-1$ terms in the summation $\Sigma''$, this shows that $\Sigma'' < (2^{k+1} - (k+1) -1)\alpha \tau(q)$,  
for any  choice of $w$ that has support $q$ satisfying condition $(3)$. Thus, if we choose $\alpha \leq 2^{-n-2}\delta$, then by  $(3)$ we get 
$\Sigma'' \leq \delta\tau(q)/2$.

So we are left with estimating the terms in the summation $\Sigma'$, which have just one occurrence of $w^j$, $j\neq 0$, 
i.e are of the form $|\tau(y_1w^jy_2)|=|\tau(w^jE_{\text{\bf A}}(qy_2y_1q))|$, for some $y_1, y_2 \in M$, 
$1\leq |j| \leq n$. There are $k+1$ many such terms for each $k= 1, ..., n$. 
Let's denote by $Y_0$ the set of all $y_1, y_2$ which appear this way, and note that this is a finite set in $q\text{\bf M}q$. Thus  $Y=E_{\text{\bf A}}(qY_0\cdot Y_0q)$ 
is finite as well. 

It is sufficient now  to find a Haar unitary $w\in \text{\bf A}_0q$ such that $|\tau(w^jy)|\leq \delta \tau(q)/2(n+1)$, $\forall y\in Y$, $1\leq |j|\leq n$, because 
then the sum of the $k+1$ terms in $\Sigma'$ will be majorized by $\delta \tau(q)/2$, altogether showing that for all $ x\in F_{u,n}$, we have $|\tau(x)|\leq \delta \tau(uu^*)$. 
Since $\text{\bf A}_0q$ is diffuse, it contains a separable diffuse von Neumann subalgebra $A_0$, which is isomorphic to $L^\infty(\Bbb T)$ with the Lebesgue 
measure corresponding to $\tau(q)^{-1}\tau_{|A_0}$. Let then $w_0\in A_0$ be a Haar unitary generating $A_0$. Since $\{w_0^m\}_m$ tends to $0$ 
in the weak operator topology and $Y\subset \text{\bf A}q$ is a finite set, there exists $n_0\geq n$ such that $|\tau(w_0^m y)|\leq \delta \tau(q)/2(n+1)$, 
for all $y\in Y$ and $|m| \geq n_0$. But then $w=w_0^{n_0}$ is still a Haar unitary and it satisfies all the required conditions. 

This ends the proof of the {\it Fact}. 

\vskip .05in 

By using this {\it Fact}, it follows that for each $n$ there exists a unitary element $v_n\in \text{\bf A}_0$ 
such that 
$$
|\tau(x)| < 2^{-n}, \forall x\in F_{v_n,n}. \tag 4
$$

For each $g\in \Gamma$, let $u_g=(u_{g,m})_m$ be a representation of $u_g$ with $u_{g,m} \in \Cal N_{M_m}(A_m)$. 
Let also $b_i=(b_{i,m})_m$ and $v_n=(v_{n,m})_m \in \text{\bf A}_0$, with $b_{i,m}, v_{n,m}\in A_m$, $\forall m$. Then $(4)$ becomes

$$
\underset m \rightarrow \omega \to \lim |\tau(b_{i,m} \Pi_{j=0}^k u_{g_{m_j},m}v_{n,m}^{n_j} u_{g_{m_j},m}^*)| < 2^{-n} \tag 5
$$
for all $1\leq i, k \leq n$, $0\leq m_0 < m_1 ... < m_k \leq n$. Also,  the fact that $v_n$ lies in $\text{\bf A}_0$ translates into  

$$
\underset m \rightarrow \omega \to \lim \| [u_{h,m}, v_{n,m}]\|_1 = 0, \forall h\in H, n \geq 1 \tag 6
$$

Let then $V_n$ be the set of all $m\in \Bbb N$ satisfying the following properties: 
 $$
  |\tau(b_{i,m} \Pi_{j=0}^k u_{g_{m_j},m}v_{n,m}^{n_j} u_{g_{m_j},m}^*)| < 2^{-n}  \tag 7 
  $$
  $$
  \| [u_{h_i,m}, v_{n,m}]\|_1 < 2^{-n}
  $$
  for all $1\leq i, k \leq n$, $0\leq m_0 < m_1 ... < m_k \leq n$, where $\{h_i\}_i=H$ is an enumeration of $H$. Note that by $(5)$ and $(6)$,  $V_n$ 
  corresponds to an open-closed neighborhood of $\omega$ in $\Omega$, under the identification $\ell^\infty(\Bbb N)=C(\Omega)$. 
  Define now recursively $W_0=\Bbb N$ and $W_{n+1}=W_n \cap V_{n+1} \cap \{n\in \Bbb N \mid n > \min W_n\}$. Then $W_n$ 
  is a strictly decreasing sequence of neighborhoods of $\omega$ (under the same identification as before) with $W_n \subset \cap_{j\leq n} V_j$.

We now define $w=(w_m)_m$, by letting $w_m=v_{n,m}$ if $m\in W_{n-1} \setminus W_n$. By the way $v_{n,m}$ have been taken,     
 $w$ follows unitary element in $\text{\bf A}$, while by the second relation in $(7)$ we have 
 $w\in \text{\bf A}^H=\text{\bf A}_0$. Also, by the first relation  in $(7)$ 
it follows that $B$ and 
$u_{g_i} \{w^n \mid n\in \Bbb Z\} u_{g_i}^*$, $i=0, 1, 2, ...$, are all multi-independent. Thus, if we denote by $A\subset \text{\bf A}$ 
the von Neumann algebra generated by $u_{g_i} \{w^n \mid n \in \Bbb Z\}u_{g_i}^*$, $i\geq 0$, then $A$ and $B$ are $\tau$-independent 
and $\Gamma \curvearrowright A$ is isomorphic to the generalized Bernoulli action $\Gamma \curvearrowright L^\infty([0,1]^{\Gamma/H})$, as desired.

\hfill 
$\square$

\proclaim{5.2. Corollary} As in $5.1$, let $A_n \subset M_n$ be a sequence of MASAs in finite factors, with $\text{\rm dim} M_n \rightarrow \infty$, and  
denote $\text{\bf A}=\Pi_\omega A_n \subset \Pi_\omega M_n = \text{\bf M}$. Let $G\curvearrowright X$ be a measure preserving 
$($but not necessarily free$)$ action of a countable amenable group $G$ on a probability space $(X,\mu)$. Let $\rho_i: L^\infty(X)\rtimes G \hookrightarrow \text{\bf M}$ 
be trace preserving embeddings taking $L^\infty(X)$ into $\text{\bf A}$, with commuting squares,  
and $G$ in the normalizer $\Cal N$ of $\text{\bf A}$ in $\text{\bf M}$, 
such that $\rho_i(G)$ acts freely on $\text{\bf A}$, $i=1,2$. Then there exists $u\in \Cal N$ such that $u\rho_1(x)u^*=\rho_2(x)$, $\forall x\in L^\infty(X)\rtimes G$. 
In particular, any two embeddings of $G$ into $\Cal N$ acting freely on $\text{\bf A}$, are conjugate by a unitary in $\Cal N$. 
\endproclaim 
\noindent
{\it Proof}. By Theorem 5.1 applied to $\Gamma=G$ and $H=\{1\}$, each one of the embeddings $\rho_i$ can be extended to embeddings, 
still denoted by $\rho_i$, of $A=L^\infty(X\times [0,1]^G)\subset L^\infty(X\times [0,1]^G) \rtimes G=M$ into $\text{\bf A}\subset \text{\bf M}$, 
satisfying the same properties, where $G \curvearrowright X \times [0,1]^G$ is the product action. This action 
is free, so the corresponding inclusion is Cartan, with $M$ AFD. Thus, by observation 2.3.2$^\circ$, the specific isomorphism $ \rho_2\circ \rho_1^{-1}: \rho_1(M) \simeq \rho_2(M)$ 
is implemented by a unitary in $\Cal N$.  
\hfill 
$\square$ 

\vskip .05in
Finally, let us mention that a slight adaption of the proof of 4.3 allows showing that given any two countable groups 
$\Gamma_1, \Gamma_2$ normalizing $D^\omega$ in $R^\omega$ (where as before $D\subset R$ is the Cartan subalgebra of the hyperfinite II$_1$ factor), 
there exists a unitary element $u\in \Cal N_{R^\omega}(D^\omega)$ that conjugates $\Gamma_1$ in free position with $\Gamma_2$. 
Moreover, if $H\subset \Gamma_1\cap \Gamma_2$ is a common amenable group, then $u$ can be taken so that to commute with $H$ 
and so that the group $\Gamma$ generated by $u\Gamma_1u^*$ and $\Gamma_2$ satisfies $\Gamma\simeq \Gamma_1 *_H \Gamma_2$,  
with $\Gamma$ acting freely if $\Gamma_1, \Gamma_2$ act freely. This recovers a result from [Pa], [ES]. We'll actually state 
and prove only the case $\Gamma_i$ act freely of such a statement, for clarity: 

\proclaim{5.3. Theorem} Let $A_n \subset M_n$ be a sequence of Cartan MASAs in finite factors, with $\text{\rm dim} M_n \rightarrow \infty$, and  
denote $\text{\bf A}=\Pi_\omega A_n \subset \Pi_\omega M_n = \text{\bf M}$, as before. Assume $\Gamma_i \subset \Cal N_{\text{\bf M}}(\text{\bf A})$ 
are  countable groups of unitaries acting freely on $\text{\bf A}$, with amenable subgroups $H_i\subset \Gamma_i$, $i=1, 2$, 
such that $H_1\simeq H_2$. Then there exists a unitary element $u\in \Cal N_{\text{\bf M}}(\text{\bf A})$ such that $uH_1u^*=H_2$ and   
such that the group generated by $u\Gamma_1u^*$ and $\Gamma_2$ is isomorphic to $\Gamma_1*_H \Gamma_2$ and acts freely on $\text{\bf A}$, where $H$ 
is the identification $H_1\simeq H_2$ under $\text{\rm Ad}(u)$.  
\endproclaim 
\noindent
{\it Proof}. By 5.2 above, there exists a unitary element $u_0\in \Cal N:=\Cal N_{\text{\bf M}}(\text{\bf A})$ 
such that $u_0H_1u_0^*=H_2$. We may thus assume $H_1=H_2$, a common subgroup we will denote by $H$. 

Denote $\text{\bf A}_0=H'\cap \text{\bf A}$.  By Corollary 5.2 there exists an $H$-invariant separable von Neumann subalgebra 
$D_0\subset \text{\bf A}$ such that the action of $H$ on $D_0$ is a Bernoulli $H$-action.  Denote $R_0=D_0\rtimes H$ and notice 
that $R_0$ is a separable, amenable von Neumann algebra, $D_0$ is a Cartan subalgebra in $R_0$, and $D_0\subset R_0$ is 
Cartan embedded into $\text{\bf A} \subset \text{\bf M}$. 
Denote $\Cal N_0=\Cal N_{R_0}(D_0)'\cap \Cal N$.

By 2.2.1$^\circ$, we have $\text{\bf A}_0'\cap \text{\bf M}=\text{\bf A}\vee H=\text{\bf A}\vee R_0$, while by 2.2.2$^\circ$, 2.2.3$^\circ$, 
we see that $\Cal N_0'\cap \Cal N=\Cal N_{R_0}(D_0)$, $\text{\bf A}_0$ is maximal abelian in $R_0'\cap \text{\bf M}$ and $\Cal N_0$  
coincides with the normalizer of $\text{\bf A}_0$ in $R_0'\cap \text{\bf M}$.

With this mind, the proof is very similar to the proof of Theorem 4.3. We will only show what the analogue of Lemma 3.2 
becomes, and leave all other details for the reader to complete. 

Denote 
$\Cal G_0=\{up \mid u\in \Cal N_0, p\in \Cal P(\text{\bf A}_0)\}$, which by 2.2.3$^\circ$ and [Dy] coincides with the set of partial isometries in $\text{\bf M}_0=\text{\bf A}_0\vee \Cal N_0$ 
that normalize $\text{\bf A}_0$.  For each finite 
subset $F \subset \Gamma_1 \cup \Gamma_2 \setminus \{1\}$, $n\geq 1$, a non-zero projection $f\in \text{\bf A}_0$ and $v\in \Cal G_0$  
satisfying $vv^*=v^*v\leq f$, we denote by   
$F_{v,n}$ the set of all elements of the form $x=u_0\Pi_{i=1}^k v^{\gamma_i}u_i$, where $u_0 \in F\cup \{1\}$, $u_i \in F$, $\gamma_i =\pm 1$, $1\leq k \leq n$. 
We need to prove  that given any $\varepsilon > 0$,  there exists  
$u\in \Cal G_0$ such that $uu^*=u^*u$, $\|E_{\text{\bf A}}(x)\|_1 \leq \varepsilon$, $\forall x\in F_{u,n}$, and $\tau(uu^*)>\tau(f)/4$. 

To do this,  let  $\delta = 2^{-n^2-1}\varepsilon$ and denote $\varepsilon_0 = \delta, \
\varepsilon_k = 2^{k+1} \varepsilon_{k-1}, \ k \geq 1$. Note that $\varepsilon_n < \varepsilon$. 
Let  $\mycal W$ denote the set of partial isometries  $v\in \Cal G_0$ with $vv^*=v^*v\leq f$ such that 
$\|E_{\text{\bf A}}(x)\|_1 \leq \varepsilon_k \tau(vv^*)$, $\forall x\in F_{v,k}$, for all $1\leq k \leq n$, and  endow $\mycal W$ with the order given by 
$w_1\leq w_2$ if $w_1=w_2w_1^*w_1$. Noticing that $\mycal W$ is well ordered with respect to $\leq $, we   
let $v\in \mycal W$ be a maximal element. Assume that $\tau(vv^*)\leq \tau(f)/4$ and note that $p=f-vv^* \in \Cal P(\text{\bf A}_0)$ will then 
satisfy $\tau(vv^*)/\tau(p) \leq 1/3$.  

If $w\in \Cal G_0$ satisfies $ww^*=w^*w\leq p$, then $u=v+w$ belongs to $\Cal G_0$ and satisfies $uu^*=u^*u$. When 
we develop $u_0\Pi_{i=1}^k(v+w)^{\gamma_i}u_i$ binomially, we get 
$$
\|E_{\text{\bf A}}(u_0\Pi_{i=1}^k u^{\gamma_i}u_i)\|_1  \leq \|E_{\text{\bf A}}(u_0 \Pi_{i=1}^k v^{\gamma_i} u_i)\|_1  + \Sigma' + \Sigma'',  
$$
where $\Sigma''$ is the sum of  the$L^1$-norm of terms that contain  at least two occurrences of $w^{\pm 1}$, while $\Sigma'$ is the sum 
the $L^1$-norm of terms containing exactly one occurrence of $w^{\pm 1}$ (as before, we use the notation $w^{-1}$ for $w^*$). 

To estimate $\Sigma''$ we use $1.4.1^\circ$, exactly the same 
way $1.4.2^\circ$ is used in the estimates $(2)-(9)$ in the proof of 3.2, to get that $\Sigma''  \leq \varepsilon_k\tau(q)-(2k+4)(\varepsilon_{k-1}/3)\tau(q)$, 
where $q=ww^*$. 
Note that in order to do that, we only use the properties of the support $q$ of $w$, namely the fact that given any 
finite set $Y\subset \text{\bf M}\ominus (\text{\bf A}\vee H)$  and any $\alpha > 0$, 
one can take $q\in \Cal P(\text{\bf A}_0)$ such that $\|qyq\|_1< \alpha \tau(q)$, $\forall y\in Y$ (by applying $1.4.1^\circ$ to $Q=\text{\bf A}_0$ 
and using the fact that $\text{\bf A}_0'\cap \text{\bf M} =\text{\bf A}\vee H$).  

Now, in order to estimate $\Sigma'$, we denote by $\Cal U_q$ the set of partial isometries in $\Cal G_0$ that have left and right support equal to 
$q$, which we view as a subgroup of unitaries in $q\text{\bf M}_0q$. Notice that $\Cal U_q$ generate $q\text{\bf M}_0q$ (by 2.2.3$^\circ$) and that $\text{\bf M}_0 
\not\prec_{\text{\bf M}} \text{\bf M}_0'\cap \text{\bf M}=R_0$ (because this centralizer is separable and amenable, and by applying 2.1 and 4.2). Thus, 
given any finite set $Y\subset \text{\bf M}$ and any $\alpha > 0$, 
there exists by 1.5 unitary elements $w\in \Cal U_q$ such that $\|E_{\text{\bf A}}(y_1wy_2)\|_1 < \alpha \tau(q)$, 
$\forall y_1, y_2\in Y$. 

Then the same estimates as the ones in $(10)-(14)$  in the proof of $3.2$, show that $u=v+w\in \mycal W$, contradicting the maximality of $v$. 
Thus, we do have indeed $\tau(vv^*) > \tau(f)/4$. With this technical fact in hand, the rest of the proof  
proceeds exactly as the proof of 4.3 in Section 4.

\hfill 
$\square$

\head  References \endhead

\item{[Be]} B. Bekka: {\it Operator-algebraic superrigidity for $SL(n, \Bbb Z)$, $n\geq 3$}, Invent. Math. {\bf 169} (2007), 401-425. 

\item{[BiJ]} D. Bisch, J.F.R. Jones: {\it Algebras associated with intermediate subfactors}, Invent. Math. {\bf 128} (1997), 89-157. 

\item{[B]} N. Brown: {\it Topological dynamical systems assciated with} II$_1$ {\it factors}, Adv. in Math. {\bf 227} (2011), 1665-1699. 

\item{[BDJ]} N. Brown, K. Dykema, K. Jung: {\it Free entropy dimension in amalgamated free products}, Proc. London Math. Soc. {\bf 97} (2008), 339-367.  

\item{[C1]} A. Connes: {\it Classification of injective factors}, Ann. of Math., {\bf 104} (1976), 73-115.

\item{[C2]} A. Connes: {\it A type} II$_1$ {\it factor with countable fundamental group}, 
J. Operator Theory {\bf 4} (1980), 151-153. 

\item{[CFW]} A. Connes, J. Feldman, B. Weiss: {\it An amenable equivalence relation is generated by a single
transformation}, Erg. Theory Dyn. Sys.  {\bf 1} (1981), 431-450.

\item{[D1]} J. Dixmier: {\it Sous-anneaux abeliens maximaux dans les facteurs de type fini}, Ann. of Math. {\bf 59} (1954), 279-286. 

\item{[D2]} J. Dixmier: ``Les alg\`ebres d'op\'erateurs dans l'espace hilbertien'', Gauthier-Vill-\newline ars, Paris 1957, 1969. 

\item{[Dy]} H. Dye: {\it On groups of measure preserving
transformations} I, Amer. J. Math, {\bf 81} (1959), 119-159.

\item{[EL]} G. Elek, G. Lippner: {\it Sofic equivalence relations}, J. Funct. Anal. {\bf 258} (2010), 1692-1708. 

\item{[ES]} G. Elek, E. Szabo: {\it  Sofic representations of amenable groups}, Proc. Amer. Math. Soc. {\bf 139} (2011), 4285-4291. 

\item{[FGR]} I. Farah, I. Goldbring, B. Hart: {\it Amalgamating $R^\omega$-embeddable von Neumann algebras}, arXiv:1301.6816. 

\item{[FHS]} I. Farah, B. Hart, D. Sherman: {\it Model theory of operator algebras} III: {\it elementary equivalence and $\text{\rm II}_1$ factors}, 
arXiv:1111.0998

\item{[F]} J. Feldman: {\it Nonseparability of certain finite factors},  Proc. Amer. Math. Soc. {\bf 7} (1956), 23-26. 

\item{[FM]} J. Feldman, C.C. Moore: {\it Ergodic equivalence
relations, cohomology, and von Neumann algebras} I, II, Trans. Amer.
Math. Soc. {\bf 234} (1977), 289-324, 325-359.

\item{[IPP]} A. Ioana, J. Peterson, S. Popa:
{\it Amalgamated Free Products of w-Rigid Factors and Calculation of their
Symmetry Groups},
Acta Math. {\bf 200} (2008), 85-153 

\item{[J]} K. Jung: {\it Amenability, tubularity, and embeddings into $R^\omega$}, Math. Ann. {\bf 338} (2007), 241-248. 

\item{[McD]} D. McDuff: {\it Central sequences and the hyperfinite factor}, Proc. London Math. Soc. {\bf 21} (1970), 443-461.

\item{[MvN1]} F. Murray, J. von Neumann:
{\it On rings of operators}, Ann. Math. {\bf 37} (1936), 116-229.

\item{[MvN2]} F. Murray, J. von Neumann: {\it Rings of operators
IV}, Ann. Math. {\bf 44} (1943), 716-808.

\item{[OW]} D. Ornstein, B. Weiss: {\it Ergodic theory of amenable group actions}, 
Bull. Amer. Math. Soc. {\bf 2} (1980), 161-164. 

\item{[O]} N. Ozawa: {\it  About the Connes embedding conjecture: algebraic approaches}, \newline arXiv:1212.1700

\item{[Pa]} L. Paunescu: {\it On sofic actions and equivalence relations}, J. Funct. Anal. {\bf 261} (2011), 2461-2485. 

\item{[Pe]} V. Pestov: {\it Hyperlinear and sofic groups: a brief guide}, Bull. Symbolic Logic {\bf 14} (2008), no. 4, 449-480. 

\item{[P1]} S. Popa: {\it On a problem of R.V. Kadison on maximal
abelian *-subalgebras in factors}, Invent. Math., {\bf 65} (1981),
269-281.

\item{[P2]} S. Popa: {\it Maximal injective subalgebras in factors
associated with free groups}, Advances in Math., {\bf 50} (1983),
27-48.

\item{[P3]} S. Popa: {\it The commutant modulo the set of compact operators of a von Neumann
algebra}, J. Funct. Analysis, {\bf 712} (1987), 393-408.

\item{[P4]} S. Popa: {\it Markov traces on universal Jones algebras and subfactors of finite index}, Invent. 
Math. {\bf 111} (1993), 375-405. 

\item{[P5]} S. Popa: {\it Classification of amenable subfactors of
type} II, Acta Mathematica, {\bf 172} (1994), 163-255.

\item{[P6]} S. Popa: {\it Free independent sequences in type} II$_1$ {\it factors
and related problems}, Asterisque, {\bf 232} (1995), 187-202.

\item{[P7]} S. Popa: {\it An axiomatization of the lattice of
higher relative commutants of a subfactor}, Invent. Math., {\bf
120} (1995), 427-445.

\item{[P8]} S. Popa: {\it The relative Dixmier property for inclusions
of von Neumann algebras of finite index}, Ann. Sci. Ec. Norm. Sup.
{\bf 32} (1999), 743-767.

\item{[P9]} S. Popa: {\it Universal construction of subfactors}, J.
reine angew. Math., {\bf 543} (2002), 39-81.

\item{[P10]} S. Popa: {\it Strong Rigidity of} II$_1$ {\it Factors
Arising from Malleable Actions of $w$-Rigid Groups} I, Invent. Math.,
{\bf 165} (2006), 369-408.

\item{[P11]} S. Popa: {\it On the superrigidity of malleable actions with spectral gap}, J. Amer. Math. Soc. {\bf 21} (2008), 981-1000. 

\item{[P12]} S. Popa: {\it A} II$_1$ {\it factor approach to the Kadison-Singer problem}, 
arXiv:1303.1424, to appear in Comm. Math. Phys. 

\item{[Va]} S. Vaes: {\it Factors of type} II$_1$ {\it without non-trivial finite index subfactors}, 
Trans. Amer. Math. Soc. {\bf 361} (2009), 2587-2606.

\item{[V]} D. Voiculescu:  {\it Symmetries of some reduced free product 
$C^*$-algebras},  In: ``Operator algebras and their connections with topology
and ergodic theory'', Lect. Notes in Math. Vol. {\bf 1132}, 566-588 (1985).

\item{[W]} B. Weiss: {\it Sofic groups and dynamical systems}, 
in ``Ergodic theory and harmonic analysis'', Mumbai, 1999, Sankhya Ser. A. {\bf 62} (2000), 350-359.

\item{[Wr]} F. B. Wright:  {\it A reduction for algebras of finite type},  Ann. of Math. {\bf 60} (1954), 560 - 570.

\enddocument